\documentclass[11pt]{amsart}

\usepackage{mathrsfs}
\usepackage{amsmath, amsthm, amsfonts}
\usepackage[percent]{overpic}

\usepackage{color}
\usepackage{MnSymbol}
\usepackage{fullpage}
\usepackage{bbm}

\newtheorem{thm}{Theorem}[section]
\newtheorem*{thm*}{Theorem}

\newtheorem{prop}[thm]{Proposition}

\theoremstyle{definition}
\newtheorem{defn}[thm]{Definition}

\newtheorem{ex}[thm]{Example}
\theoremstyle{remark}
\newtheorem{rem}[thm]{Remark}

\newcommand{\QQ}{\mathbb{Q}}
\newcommand{\EE}{\mathbb{E}}

\newcommand{\HH}{\mathcal{H}}
\newcommand{\PPP}{\mathcal{P}}

\newcommand{\RR}{\mathbb{R}}
\newcommand{\CC}{\mathbb{C}}

\newcommand{\OO}{\mathcal{O}}
\newcommand{\ord}{\operatorname{ord}}

\newcommand{\ZZ}{\mathbb{Z}}

\newcommand{\PP}{\mathbb{P}}

\newcommand{\dJ}{\text{dJ}}

\newcommand{\Pic}{\text{Pic}}
\newcommand{\Mg}{\mathcal{M}_g}
\newcommand{\Mgn}{\mathcal{M}_{g,n}}
\newcommand{\Mgb}{\overline{\mathcal{M}}_g}
\newcommand{\Mgnb}{\overline{\mathcal{M}}_{g,n}}
\newcommand{\Mgm}{\mathcal{M}_{g,1}}
\newcommand{\Mgbm}{\overline{\mathcal{M}}_{g,1}}

\newcommand{\rk}{\text{rk}}

\newcommand{\res}{\text{res}}

\newcommand{\SL}{\mbox{SL}}

\newcommand{\Div}{\operatorname{Div}}

\newtheorem{Thm*}{Theorem*}
\theoremstyle{definition}

\title{Effective divisors in $\Mgnb$ from abelian differentials}
\author{Scott Mullane}
\date{\today}

\AtEndDocument{\bigskip{\footnotesize%
\textsc{Scott Mullane, Department of Mathematics, Harvard University, One Oxford St, Cambridge, MA 02138, USA} \par
\textit{E-mail address}: \texttt{smullane@math.harvard.edu} \par
}}

\begin{document}
\thispagestyle{empty}

\maketitle

\begin{abstract}
We compute many new classes of effective divisors in $\Mgnb$ coming from the strata of abelian differentials and efficiently reproduce many known results obtained by alternate methods. Our method utilises maps between moduli spaces and the degeneration of abelian differentials.
\end{abstract}

\setcounter{tocdepth}{1}

\tableofcontents

\section{Introduction}
 The moduli space of abelian differentials $\HH(\kappa)$ consists of pairs $(C,\omega)$ where $\omega$ is a holomorphic or meromorphic differential on a smooth genus $g$ curve $C$ and the multiplicity of the zeros and poles of $\omega$ is fixed of type $\kappa$, an integer partition of $2g-2$. Previous seminal work has exposed the fundamental algebraic attributes of these spaces~\cite{KontsevichZorich}\cite{McMullen}\cite{EskinMirzakhani}\cite{EskinMirzakhaniMohammadi}. From the perspective of algebraic geometry, geometrically defined codimension one subvarieties, or divisors, have been used to study many aspects of moduli spaces of curves including the Kodaira dimension and the cone of effective divisors~\cite{HarrisMumford}\cite{EisenbudHarrisKodaira}\cite{Farkas23}
 . In this paper we compute the class of many effective divisors in $\Mgnb$ defined by the strata of abelian differentials.

The divisor $D^n_\kappa$  in $\Mgnb$ for $\kappa=(k_1,...,k_m)$ with $\sum k_i=2g-2$ is defined as 
\begin{equation*}
D^n_\kappa=\overline{\{[C,p_1,...,p_n]\in\Mgn\hspace{0.15cm}|\hspace{0.15cm} [C,p_1,...,p_m]\in\mathcal{M}_{g,m} \text{ with  }  \sum k_ip_i\sim K_C   \}},
\end{equation*}
where $m=n+g-2$ or $n+g-1$ for holomorphic and meromorphic signature $\kappa $ respectively. When all $k_i$ are even this divisor has two irreducible components based on spin structure that we denote by the indices \emph{odd} and \emph{even}. We use this notation to denote both the divisor and the class\footnote{Our class expressions are given modulo the labelling of the unmarked points. See Section~\ref{notation} for an explicit description.} of the divisor in $\Pic(\Mgnb)\otimes \QQ$.

The main results of Section~\ref{Coupled}, Section~\ref{pp} and Section~\ref{Mgbm} are summarised in the following three theorems.

\begin{thm}[Coupled partition divisors] 
\begin{equation*}
D^2_{1,1,2^{g-2}}=2^{g-3}(2^{g+1}\lambda +2^{g-1}(\psi_1+\psi_2) -2^{g-2}\delta_0    -\sum_{i=0}^{g-1}2^{i+1}(2^{g-i}-1)\delta_{i:\{1,2\}}      -\sum_{i=1}^{g-1}2^{g-1}\delta_{i:\{1\}}         );
\end{equation*}
For $\underline{d}=(d_1,...,d_{n})$ such that $\sum_{j=1}^{n}d_j=0$ with  $\underline{d}^-\ne \{-2\}$, then
\begin{equation*}
D^{n}_{\underline{d},2^{g-1}}=2^{g-2}(2^{g+1}\lambda +2^{g-1}\sum_{j=1}^nd_j^2\psi_j-2^{g-2}\delta_0   
 -\sum_{\tiny{\begin{array}{cc}|d_S|=0\\|S|\ne n  \end{array}}}\sum_{i=0}^{g}2^{g-i+1}(2^i-1)\delta_{i:S}   
 -2^{g-1}\sum_{|d_S|\geq1}\sum_{i=0}^{g-1}d_S^2\delta_{i:S}        
 ).
\end{equation*}
If all $d_j$ are even then
\begin{equation*}
D^{n,\textit{odd}}_{\underline{d},2^{g-1}}=2^{g-2}((2^g-1)\lambda+\frac{2^g-1}{4}\sum_{j=1}^n d_j^2\psi_j-2^{g-3}\delta_0
-\sum_{\tiny{\begin{array}{cc}|d_S|=0\\|S|\ne n  \end{array}}}\sum_{i=0}^{g}(2^i-1)(2^{g-i}+1)\delta_{i:S}   -\frac{2^g-1}{4}\sum_{|d_S|\geq 2}\sum_{i=0}^{g-1}d_S^2\delta_{i:S}         )
\end{equation*}
and 
\begin{equation*}
D^{n,\textit{even}}_{\underline{d},2^{g-1}}=2^{g-2}((2^g+1)\lambda +\frac{2^g+1}{4}\sum_{j=1}^n d_j^2\psi_j -2^{g-3}\delta_0
-\sum_{|d_S|=0}\sum_{i=0}^{g}(2^i-1)(2^{g-i}-1)\delta_{i:S}     -\frac{2^g+1}{4}\sum_{|d_S|\geq 2}\sum_{i=0}^{g-1}d_S^2\delta_{i:S}         ).
\end{equation*}
For $\underline{d}=(-2,1,1)$, 
\begin{eqnarray*}
D^{3}_{\underline{d},2^{g-1}}&=&2^{g-3}(2^{g+1}\lambda +2^{g+2}\psi_1+2^{g-1}(\psi_2+\psi_3)-2^{g-2}\delta_0 -\sum_{i=0}^{g-1}2^{i+1}(2^{g-i}-1)\delta_{i\{1,2,3\}}-\sum_{i=0}^{g-1}2^{i+1}(2^{g-i}+1)\delta_{i:\{2,3\}}\\
&&-\sum_{i=0}^{g-1}2^{g-1}(\delta_{i:\{1,2\}}+\delta_{i:\{1,3\}}).
\end{eqnarray*}
For $\underline{d}=(-2,2)$, then
\begin{equation*}
D^{2}_{\underline{d},2^{g-1}}=2^{g-3}(2^{g+1}\lambda+2^{g+2}\psi_1 +2(2^g+1)\psi_2-2^{g-2}\delta_0 -\sum_{i=0}^{g-1}2^{i+1}(2^{g-i}-1)\delta_{i\{1,2\}}-2^{i+1}(2^{g-i}+1)\sum_{i=1}^{g-1}\delta_{i:\{2\}} ),
\end{equation*}
with
\begin{equation*}
D^{2,\textit{odd}}_{\underline{d},2^{g-1}}=\varphi_2^*D^{1,\textit{odd}}_{2^{g-1}}=2^{g-3}((2^g-1)\lambda+2(2^g-1)\psi_1+0\psi_2-2^{g-3}\delta_0-\sum_{i=0}^{g-1}(2^i+1)(2^{g-i}-1)\delta_{i:\{1,2\}}-\sum_{i=1}^{g-1}(2^i-1)(2^{g-i}+1)\delta_{i:\{2\}})
\end{equation*}
where $\varphi_2:\overline{\mathcal{M}}_{g,2}\longrightarrow\Mgbm$ forgets the second marked point, and
\begin{equation*}
D^{2,\textit{even}}_{\underline{d},2^{g-1}}=2^{g-3}((2^g+1)\lambda+2(2^g+1)\psi_1 +2(2^g+1)\psi_2 -2^{g-3}\delta_0-\sum_{i=0}^{g-1}(2^i-1)(2^{g-i}-1)\delta_{i:\{1,2\}} -\sum_{i=1}^{g-1}(2^i+1)(2^{g-i}+1)\delta_{i:\{2\}} ).
\end{equation*}
\end{thm}

\begin{thm}[Pinch partition divisors]
For $\underline{d}=(d_1,...,d_n)$ with $d_i\geq0$ and $\sum d_i=g-1$ we have for $g\geq 3$
\begin{equation*}
D^n_{\underline{d},1^{g-3},2}=-4(g-7)\lambda+\sum_{i=1}^{g-1}(2g(d_i+1)-3d_i-5)d_i\psi_i-2\delta_0+\sum_{i=0}^g\sum_{d_S=0}^{i-1}c_{i:S}\delta_{i:S}
\end{equation*}
where\footnote{Note that $c_{0:S}=c_{g:g-d_S-1}$ for $|S|\geq 2$.}
\begin{equation*}
c_{i:S}=(3-2g)d_S^2+(4gi+2g-10i+1)d_S -2gi^2+7i^2-2gi-i-2,
\end{equation*}
for $d_S:=\sum_{i\in S}d_i$. 

For $\underline{d}=(d_1,...,d_n)$ with $\sum d_i=g-2$, $d_j\leq -2$ and $d_i\geq0$ for $i\ne j$, then for $d_j\leq -3$,
\begin{eqnarray*}
D^{n}_{\underline{d},1^{g-2},2}&=&(26-4g)\lambda +\sum_{i=1}^n 2d_i((g-1)d_i+g-2)\psi_i-2\delta_0    +\sum_{i=0}^{g-1}c_{i:S}\delta_{i:S}
\end{eqnarray*}
where for $j\nin S$ and $d_S\leq i-1$,
\begin{equation*}
c_{i:S}=(2-2g) d_S^2+2 ( 2g i+g-4i+1) d_S-2 (g i^2+g i-3i^2 +i+1)
\end{equation*}
and for $d_S\geq i$,
\begin{equation*}
c_{i:S}=(2-2g) d_S^2  +2 (2g i-g-4i+2) d_S -2 (g i^2-3i^2-g i+4i).
\end{equation*}
For $d_j=-2$
\begin{eqnarray*}
D^{n}_{\underline{d},1^{g-2},2}&=&(27-4g)\lambda +4g\psi_j+\sum_{i\ne j} \frac{(4g(d_i+1)-5d_i-9)d_i}{2}\psi_i-2\delta_0    +\sum_{i=0}^{g-1}c_{i:S}\delta_{i:S}
\end{eqnarray*}
where for $j\nin S$ and $d_S\leq i-1$,
\begin{equation*}
c_{i:S}=\frac{1}{2}((5-4g)d_S^2+(8g i+4g-18i+3)d_S-4g i^2-4gi+13i^2-3i-4)
\end{equation*}
and for $d_S\geq i$,
\begin{equation*}
c_{i:S}=\frac{1}{2}((5-4g)d_S^2+(8g i-4g-18i+9)d_S-4g i^2+4g i+13i^2-17i).
\end{equation*}
\end{thm}

\begin{thm}[Divisors in $\Mgbm$]\label{thmMgbm}
The divisor $D^1_{g-k,k+1,1^{g-3}}$ for $g\geq3$ and $k=1,...,g-1$ is given by\footnote{This formula  is reproduced in the case that $k=g-1$ in the Appendix Section~\ref{sec:appendix} by the more labour intensive methods of Porteous' formula and test curves. Setting $k=0$ in this formula recovers $(g-2)W$ where $W$ is the Weierstrass divisor computed  by Cukierman~\cite{Cukierman}. 
This divisor is irreducible in all cases except $g=3$, $k=1$ which is discussed in Remark~\ref{rem:22}.}
\begin{eqnarray*}
c_\psi&=&\frac{(k+1)(g-k)((k+1)g^2-(k^2+k+1)g-2)}{2}\\
c_\lambda&=& \frac{ (k+1) (4 - 2 g + 10 k - 2 g k + 11 k^2 + 3 k^3)}{2}\\
 c_0&=&\frac{(k+1)^2-(k+1)^4}{6}\\
 c_i&=&\begin{cases}- \frac{(k+1)(i(g-i+1)(g-i)(k+1)+(g-i-k)((k+1)(g-i)^2-(k^2+k+1)(g-i)-2))}{2}&\text{ for $1\leq i\leq g-k$.}\\
 - \frac{(g - i) (k+1) (-3 g + g^2 + 4 i - g i + 3 k - 4 g k + g^2 k + 
   5 i k - g i k + 3 k^2 - 2 g k^2 + 2 i k^2 + k^3)}{2} &\text{ for $g-k+1\leq i\leq g-1$.}
 \end{cases}
 \end{eqnarray*}
The divisor $D^1_{-h,g+h,1^{g-2}}$ for $h\geq 2$ is given by\footnote{When $g=2$ and $h$ is even this divisor has two connected components based on spin structure. $D^{1,\textit{odd}}_{-h,h+2}=5W$ for $4$ $|$ $h$ and $D^{1,\textit{even}}_{-h,h+2}=5W$ otherwise.}
\begin{eqnarray*}
c_\psi&=&\frac{g(g+h+1)(h-1)(h^2+gh+g+1)}{2}\\
c_\lambda&=&\frac{(1 + g + h) (2 - 3 g^2 + 3 g^3 - 2 h - 4 g h + 9 g^2 h - h^2 + 
   9 g h^2 + 3 h^3)}{2}\\
 c_0&=&\frac{(g+h)^2-(g+h)^4}{6}\\
 c_i&=&\frac{(i-g)}{2} (g^3(2h+i)+g^2(5h^2+2i h+2h+i-1)+g(4h^3+h^2(i+4)+2h(i-1)+i)\\
 &&-2h-i+1+h(h+2)(i+h^2))\hspace{0.3cm}\text{ for $1\leq i\leq g-1$.}
\end{eqnarray*}
\end{thm}

In~\cite{Mullane} the author used the method of test curves to obtain a closed formula for all divisors $D^0_\kappa$ in $\Mgb$ for all holomorphic signatures $\kappa$. The setting of $\Mgnb$ allows us the opportunity to provide an exposition of a different method of calculating divisor classes.  In this paper we employ maps between moduli spaces to compute the classes of interest. Bainbridge, Chen, Gendron, Grushevsky, and M\"oller~\cite{BCGGM} have recently provided a full compactification of the space of abelian differentials.
With this understanding of the degeneration of differentials we are able to explicitly describe the components of the pullback of a divisor coming from the strata of differentials under the maps described in Section~\ref{sec:maps} obtained by gluing in marked tails of different genus at marked points and gluing marked points together. Hence by knowing the class of the components of the pullback of an unknown divisor class or identifying an unknown divisor as a component of the pullback of a known class, we obtain many of the coefficients or the full class of the unknown divisor.

The difficulties in the computation arise from obtaining sufficient relations to find all coefficients 
and from computing the multiplicity of each component of the pullback of a divisor in such a relation. Obtaining the multiplicity requires enumerating certain holomorphic and meromorpic sections of a line bundle on a general curve. We use a variant of the well-known de Jonquieres' formula in the holomorphic case. In the meromorphic case the Picard variety method realises the unknown number as the degree of a map between $C^g$ and the Picard variety. The difficulty of calculating the multiplicity of different specific solutions that violate the global residue condition can be overcome by investigating the ramification locus of this map. 

Section~\ref{Mgbm} provides an illustrative introduction to the techniques that will be developed in later sections by generalising the Weierstrass divisor $W=D^1_{g,1^{g-2}}$, the closure in $\Mgbm$ of the locus of Weierstrass points originally calculated by Cukierman~\cite{Cukierman}.
The families of divisors computed in Section~\ref{Coupled} and Section~\ref{pp} were chosen to best expose the utility of our methods and are relevant in our search for extremal effective divisor classes. The coupled partition divisors present the simplest case to provide an exposition of our method as it pertains to strata of abelian differentials with multiple components. When all the marked points have even multiplicities, these divisors have two components based on spin structure. The pinch partition divisors provide a useful exposition of our method in the case that the unmarked points have different multiplicities and the first example of a pinch partition divisor, $D^{g-1}_{1^{2g-4},1}$, was shown to be extremal in the effective cone by Farkas and Verra~\cite{FarkasVerraU}.

Divisors in $\Mgnb$ from the strata of abelian differentials have been presented previously in various places, though often under different guises. In Section~\ref{known} we study the known classes coming from strata of abelian differentials. Logan~\cite{Logan} investigated the Kodaira dimension of $\Mgnb$ through the use of effective divisors defined  for $\underline{d}=(d_1,...,d_n)$, $d_i\geq0$ with $\sum d_i=g$ as the closure of $[C,p_1,...,p_n]\in \Mgn$ such that $|\sum d_ip_i|$ is a $g^1_g$. In Section~\ref{known} we show that these are the divisors $D^n_{\underline{d},1^{g-2}}$ and in Section~\ref{ReplicatingLogan} we reproduce these calculations from this perspective.  

M\"uller~\cite{Muller}, Grushevsky and Zakharov~\cite{GruZak} computed class of the closure of the pullback of the theta divisor, that is, for $\underline{d}=(d_1,...,d_n)$ with $\sum d_i=g-1$ and some $d_i<0$, the closure of $[C,p_1,...,p_n]\in \Mgn$ such that $h^0(\sum d_ip_i)\geq1$. In Section~\ref{known} we show that these are the divisors $D^n_{\underline{d},1^{g-1}}$ generalising Logan's divisors to the meromorphic case. In Section~\ref{ReplicatingMuller} we reproduce these calculations from this perspective. 

The class of strata with more than one irreducible connected component has also been considered in one isolated case. Teixidor i Bigas~\cite{Teixidor} computed the closure of the locus of curves in $\Mgb$ with a vanishing theta null or curves that admit a semi-canonical pencil. Pulling this locus back under the map $\varphi:\Mgbm\longrightarrow \Mgb$ we obtain the locus of points on a curve with a vanishing theta null or the even component of the divisor $D^1_{2^{g-1}}$ which we denote $D^{1,\textit{even}}_{2^{g-1}}$. Farkas and Verra~\cite{FarkasVerraTheta} computed the odd spin structure component $D^{1,\textit{odd}}_{2^{g-1}}$ as the closure of the loci of points in the support of an odd theta characteristic. In Section~\ref{ReplicatingFarkasTexidor} we show how each of these results implies the other and how to simply reproduce the class calculation through the use of maps between moduli spaces and the degeneration of differentials and theta characteristics to nodal curves.

Farkas and Verra~\cite{FarkasVerraTheta} computed the class of the closure of the anti-ramification locus in $\overline{\mathcal{M}}_{g,g-1}$ for $g\geq 3$ defined as the closure of $[C,p_1,...,p_{g-1}]\in \mathcal{M}_{g,g-1}$ such that $h^0(K_C-p_1-...-p_{g-1}-2q)\geq 1$ for some $q\in C$. From our perspective this divisors is $D^{g-1}_{1^{2g-4},2}$. In Section~\ref{ReplicatingFarkasVerra} we efficiently replicate this class calculation by our alternate method. 
\\
\\
\textbf{Acknowledgements.} This project formed part of my PhD thesis. I would like to thank my advisor Dawei Chen for his guidance and many helpful discussions. During the preparation of this article, the author was partially supported by CAREER grant 1350396.

\section{Preliminaries}\label{Prelims}

\subsection{Strata of abelian differentials}\label{strata:abelian}
A partition of $2g-2$ of the form $\kappa=(k_1,...,k_n)$ with all $k_i\in \ZZ\setminus\{0\}$ is known as a \emph{signature}. We say $\kappa$ is \emph{holomorphic}, denoted $\kappa>0$ if all entries $k_i>0$ and $\kappa$ is \emph{meromorphic}, denoted $\kappa\not>0$ if some $k_i<0$. We define the \emph{stratum of abelian differentials with signature $\kappa$} as
\begin{equation*}
\HH(\kappa):=\{ (C,\omega) \hspace{0.15cm}|\hspace{0.15cm} g(C)=g,\hspace{0.05cm} (\omega)=k_1p_1+...+k_np_n, \text{ for $p_i$ distinct}\}
\end{equation*}
where $\omega$ is a meromorphic differential on $C$. Hence $\HH(\kappa)$ is the space of abelian differentials with prescribed multiplicities of zeros and poles given by $\kappa$. By relative period coordinates $\HH(\kappa)$ has dimension $2g+n-1$ if $\kappa>0$ and $2g+n-2$ if $\kappa\not>0$. 

A related object of interest in our study of the birational geometry of $\Mgn$ is the \emph{stratum of canonical divisors with signature $\kappa$} which we define as
\begin{equation*}
\PPP(\kappa):=\{[C,p_1,...,p_n]\in \Mgn   \hspace{0.15cm}| \hspace{0.15cm}k_1p_1+...+k_np_n\sim K_C\}.
\end{equation*}
Forgetting this ordering of the zeros or poles of the same multiplicity we obtain the projectivisation of $\HH(\kappa)$. If all $k_i$ are distinct then this finite cover becomes an isomorphism of $\PPP(\kappa)$ with the projectivisation of $\HH(\kappa)$ under the $\CC^*$ action that scales the differential $\omega$. Hence we have the dimension of $\PPP(\kappa)$ is $2g+n-2$ if $\kappa>0$ and $2g+n-3$ if $\kappa\not>0$. 

A line bundle $\eta$ on a smooth curve $C$ such that $\eta^{\otimes 2}\sim K_C$ is known as a \emph{theta characteristic}. The \emph{spin structure} of $\eta$ is the parity of $h^0(C,\eta)$ which Mumford~\cite{Mumford} showed is deformation invariant. Consider an abelian differential $(C,\omega)$ where $\omega$ has signature $\kappa=(k_1,..,k_n)$. If all $k_i$ are even then an abelian differential of this type specifies a theta characteristic on the underlying curve
\begin{equation*}
\eta\sim \sum_{i=1}^n \frac{k_i}{2}p_i.
\end{equation*}
As the parity of $h^0(C,\eta)$ is deformation invariant, the loci $\HH(\kappa)$ and $\PPP(\kappa)$ are reducible and break up into disjoint components with even and odd parity of $h^0(C,\eta)$. 

A signature $\kappa>0$ is of \emph{hyperelliptic type} if all odd entires in the signature occur in pairs $\{j,j\}$ or $\{-j,-j\}$. A \emph{hyperelliptic differential} of type $\kappa$ for such a $\kappa$, is a differential on a hyperelliptic curve resulting from pulling back a degree $g-1$ rational function under the unique hyperelliptic cover of $\PP^1$, with the minimum number of zeros occurring at ramification points of the hyperelliptic involution known as \emph{Weierstrass points}. Hence the subvariety of hyperelliptic differentials in $\HH(\kappa)$ has dimension $2g+(n-m)/2$ where $m$ is the number of zeros that occur at Weierstrass points in each hyperelliptic differential and this is the subvariety of maximum dimension that can be built in $\HH(\kappa)$ from the locus of hyperelliptic curves. Kontsevich and Zorich~\cite{KontsevichZorich} showed that there can be at most $3$ connected components in total of $\HH(\kappa)$ for $\kappa>0$ and hence $\PPP(\kappa)$, corresponding to the case that the hyperelliptic differentials become a connected component of $\HH(\kappa)$ distinct from the remaining differentials that provide two further connected components based on odd or even spin structure. Boissy~\cite{Boissy} showed that this holds in the meromorphic case ($\kappa\not>0$) for $g\geq 2$ and completely classified the connected components of $\HH(\kappa)$ when $g=1$.

\subsection{Degeneration of abelian differentials}
The investigation of how abelian differentials degenerate as the underlying curve becomes singular has recently attracted much attention. In calculating the Kodaira dimension of a number of the strata of abelian differentials, Gendron~\cite{Gendron} used analytic methods to investigate the degeneration of abelian differentials. Chen~\cite{ChenDegen} used algebraic methods to consider the limiting position of Weierstrass points on general curves of compact type. Farkas and Pandharipande~\cite{FP} extended these ideas to all nodal curves defining \emph{the moduli space of twisted canonical divisors of type $\kappa$} and showed that this space was, in general, reducible and contained extra boundary components.  Janda, Pandharipande, Pixton and Zvonkine in an appendix to this paper provided a conjectural description of the cohomology classes of the strata. Concurrently in~\cite{Mullane}, the author obtained a closed formula for the class of the strata closure that form a codimension one subvariety in $\Mgb$.

A twisted canonical divisor of type $\kappa=(k_1,...,k_n)$ is a collection of (possibly meromorphic) divisors $D_j$ on each irreducible component $C_j$ of $C$ such that
\begin{enumerate}
\item
The support of $D_j$ is contained in the set of marked points and the nodes lying in $C_j$, moreover if $p_i\in C_j$ then $\ord_{p_i}(D_j)=k_i$.
\item
If $q$ is a node of $C$ and $q\in C_i\cap C_j$ then $\ord_q(D_i)+\ord_q(D_j)=-2$.
\item
If $q$ is a node of $C$ and $q\in C_i\cap C_j$ such that $\ord_q(D_i)=\ord_q(D_j)=-1$ then for any $q'\in C_i\cap C_j$ we have $\ord_{q'}(D_i)=\ord_{q'}(D_j)=-1$. We write $C_i\sim C_j$.
\item
If $q$ is a node of $C$ and $q\in C_i\cap C_j$ such that $\ord_q(D_i)>\ord_q(D_j)$ then for any $q'\in C_i\cap C_j$ we have $\ord_{q'}(D_i)>\ord_{q'}(D_j)$. We write $C_i\succ C_j$.
\item
There does not exist a directed loop $C_1\succeq C_2\succeq...\succeq C_k\succeq C_1$ unless all $\succeq$ are $\sim$.
\end{enumerate}  
The natural question is what other conditions are required to distinguish the main component coming from twisted canonical divisors on smooth curves from the boundary components. Bainbridge, Chen, Gendron, Grushevsky, and M\"oller~\cite{BCGGM} have recently provided the global residue condition required to distinguish the main component from the boundary components giving a full compactification for the strata of abelian differentials. Let $\Gamma$ be the dual graph of $C$. They show that a twisted canonical divisor of type $\kappa$ is the limit of twisted canonical divisors on smooth curves if there exists a collection of meromorphic differentials $\omega_i$ on $C_i$ with $\Div(\omega_i)=D_i$ that satisfy the following conditions
\begin{enumerate}
\item
If $q$ is a node of $C$ and $q\in C_i\cap C_j$ such that $\ord_q(D_i)=\ord_q(D_j)=-1$ then $\res_q(\omega_i)=\res_q(\omega_j)$.
\item
There exists a full order on the dual graph $\Gamma$, written as a level graph $\overline{\Gamma}$, agreeing with the order of $\sim$ and $\succ$, such that for any level $L$ and any component $Y$ of  $\overline{\Gamma}_{>L}$ that does not contain a prescribed pole we have
\begin{equation*}
\sum_{\begin{array}{cc}\text{level}(q)=L, \\q\in C_i\in Y\end{array}}\res_{q}(\omega_i)=0
\end{equation*} 
\end{enumerate}
Part (b) is known as the \emph{global residue condition (GRC)}. Consider the following example.

\begin{ex}
Consider the nodal genus $g=4$ curve $C$ with three irreducible components $g(X)=g(Y)=0$ and $g(Z)=3$ and $X\cap Z=\{q_1,q_2\}$, $Y\cap Z=\{q_3\}$ and all other intersections zero.
Figure 1 depicts a twisted canonical divisor of type $\kappa=(6,-1,5,-4)$ on $C$ with
\begin{eqnarray*}
D_X&=&6p_1-p_2-4q_1-3q_2,\\
D_Y&=&5p_3-4p_4-3q_3,\\
D_Z&=&2q_1+q_2+q_3,
\end{eqnarray*}
where $D_j\sim K_{j}$. 
\begin{figure}[htbp]
\begin{center}
\begin{overpic}[width=0.6\textwidth]{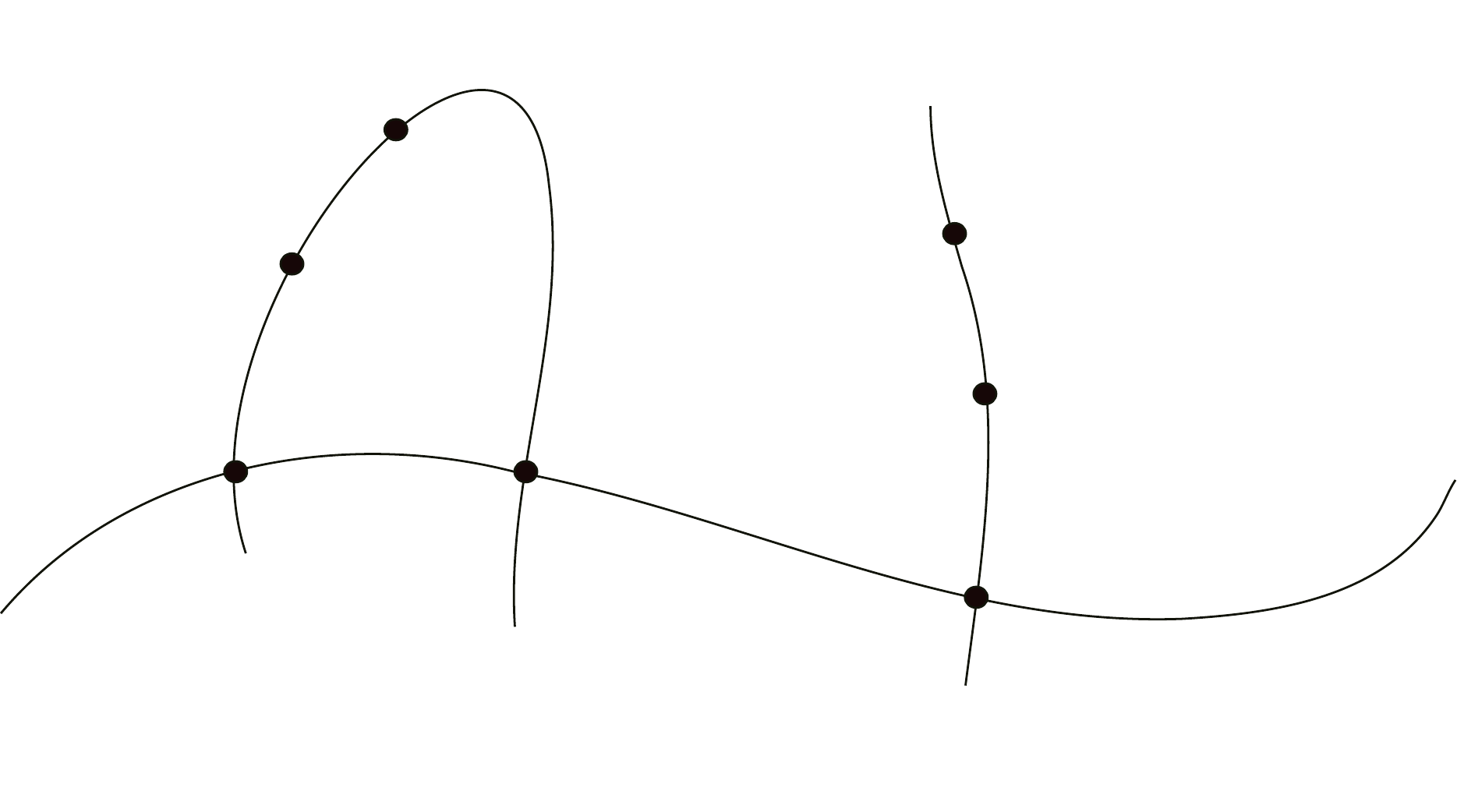}
\put (15, 38){$p_1$}
\put (22,48){$p_2$}
\put (10,24){$$}
\put (36,24){$$}
\put (12,20){$q_1$}
\put (36,20){$q_2$}
\put (63,12){$q_3$}
\put (62,17){$$}
\put (61,28){$p_4$}
\put (60,39){$p_3$}

\put (-5,16){$Z$}

\put (38,50){$X$}
\put (65,50){$Y$}

\put (-5,-0){\footnotesize{Figure 1: A twisted canonical divisor of type $\kappa=(6,-1,5,-4)$ on a nodal curve }}
\end{overpic}
\end{center}
\end{figure}
To find the conditions on such a twisted canonical divisor being smoothable we must consider the possible level graphs of the components and see what conditions these will place on the residues. 
Consider the meromorphic differential on $Y\cong \PP^1$ with signature $(5,-4,-3)$. By the cross ratio we can set the poles to $0$ and $\infty$ and the zero to $1$. The resulting differential is given locally at $0$ by
\begin{equation*}
c\frac{(1-z)^{5}}{z^{4}}dz
\end{equation*}  
for some constant $c\in \CC^*$. Hence the residue at $0$ and $\infty$ are non-zero. The flat geometric way of presenting this is that on a genus $g=0$ surface with one conical singularity, the length of any saddle connection is obtained by integrating the differential along the path which is equal to the sum of the residues the path encloses. Hence no such surface can have zero residues at both poles and we will have $\res_{q_3}(\omega_Y)\ne 0$. This shows that in the level graph the components $X$ and $Y$ must sit at the same level  and the only possible level graph is shown in Figure 2. The global residue condition on this level graph becomes
\begin{equation*}
\res_{q_1}(\omega_X)+\res_{q_2}(\omega_X)+\res_{q_3}(\omega_Y)=0.
\end{equation*}
By the residue theorem we know that the sum of the residues on any $\omega_j$ is zero and hence our condition is equivalently
\begin{equation*}
\res_{q_3}(\omega_Y)=\res_{p_2}(\omega_X).
\end{equation*}
\begin{figure}[htbp]
\begin{center}
\begin{overpic}[width=0.25\textwidth]{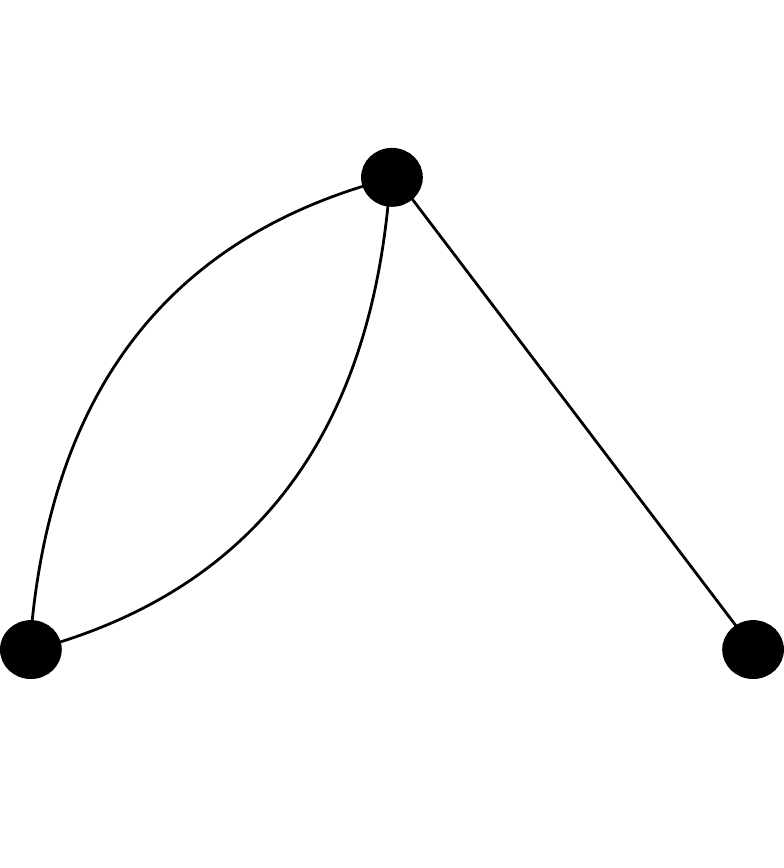}
\put (-1, 13){$X$}
\put (82,13){$Y$}
\put (39, 87){$Z$}

\put (-50,-0){\footnotesize{Figure 2: The level graph giving the global residue condition.}}
\end{overpic}
\end{center}
\end{figure}
We have seen $\res_{q_3}(\omega_Y)\ne 0$ and as $p_2$ is a simple pole, $\res_{p_2}(\omega_X)\ne 0$. Hence as there exist $\omega_i$ that satisfy $\Div(\omega_i)=D_i$, by scaling we can always satisfy this condition and we have shown that all twisted canonical divisors of this type are smoothable.

By investigating topologically a family of twisted canonical divisors on smooth curves degenerating to a nodal curve we can see why this condition on the residues is necessary. Let $\chi$ be a family of meromorphic differentials $(C_t,\omega_t)$ with $\omega_t$ of type $\kappa=(6,-1,5,-4)$ on smooth curves $C_t$ for $t\ne 0$, degenerating to the nodal curve $C_0=C$ at $t=0$. Figure $3$ depicts topologically an element of this family for $t\ne 0$. Let $v_i$ for $i=1,2,3$ be the vanishing cycles on $X_t\cup Y_t\cup Z_t$ that shrink to the nodes $q_i$ at $t=0$ such that $X_t\cap Z_t=\{v_1,v_2\}$ and $Y_t\cap Z_t=v_3$ with $X_t\to X,Y_t\to Y,$ and $ Z_t\to Z$ as $t\to 0$. As there are no poles on the component $Z_t$ we observe by an application of Stokes formula to the cycle $v_4$ that 
\begin{equation*}
\int_{v_1+v_2+v_3}\omega_t=\int_{v_4}\omega_t=0
\end{equation*}
for $t\ne 0$. Our residue condition is simply the limit of this condition as $t\to 0$.
This shows that this residue condition is necessary. Complex-analytic plumbing techniques and flat geometry are used in~\cite{BCGGM} to show in all cases the condition is sufficient.
\begin{figure}[htbp]
\begin{center}
\begin{overpic}[width=0.5\textwidth]{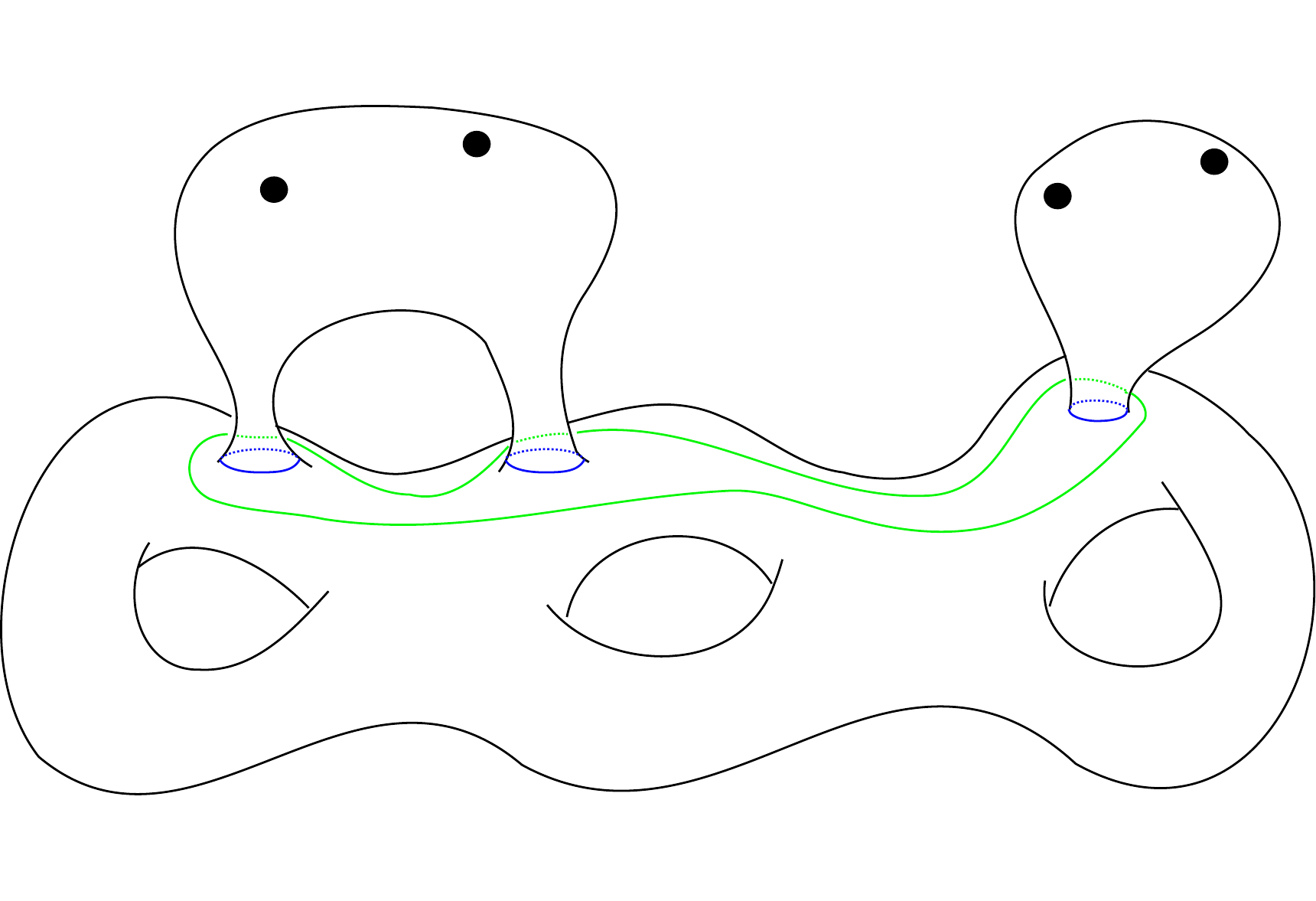}
\put (19, 29){$v_1$}
\put (41,29){$v_2$}
\put (83,33){$v_3$}
\put (67,25){$v_4$}

\put (19, 49){$p_1$}
\put (35, 52){$p_2$}
\put (80, 49){$p_3$}
\put (91, 51){$p_4$}

\put (-8,16){$Z_t$}

\put (6,50){$X_t$}
\put (70,50){$Y_t$}

\put (-15,-0){\footnotesize{Figure 3: A twisted canonical divisor of type $\kappa=(6,-1,5,-4)$ on a smooth curve }}
\end{overpic}
\end{center}
\end{figure}

\end{ex}

\subsection{Degeneration of theta characteristics and spin structures}\label{theta}
Distinguishing how different components of the strata of canonical divisors extend to the boundary of $\Mgnb$ will depend on understanding how theta characteristics and spin structures degenerate. Cornalba \cite{Cornalba} investigated how theta characteristics degenerate to nodal curves of pseudocompact type. We begin by considering a curve $C$ with a non-separating node which is the only node in the curve. Let $\tilde{C}$ be the normalisation of $C$ and $x$ and $y$ be such that $C=\tilde{C}/\{x\sim y\}$. There are two types of theta characteristic on such a curve. Consider $\tilde{\eta}$ on $\tilde{C}$ such that $\tilde{\eta}^{\otimes 2}\sim K_{\tilde{C}}+x+y$. As $K_X+x$ has a base point $x$ for any $x$, a section of $H^0(\tilde{C},\tilde{\eta})$ vanishes at $x$ if and only if it vanishes at $y$. Hence such sections are codimension one in $H^0(\tilde{C},\tilde{\eta})$. For any $\tilde{\eta}$ there are two ways to glue sections to agree at the node and hence descend to $C$. Sections can be glued as $f(x)=f(y)$ or $f(x)=-f(y)$. Hence $h^0(C,\eta)$ will differ by $+1$ for these two cases, representing even and odd spin structures. This obtains $2^{2g-2}$ even and $2^{2g-2}$ odd spin structures on the curve $C$. 

By blowing up at the node and inserting a $\PP^1$-bridge between $x$ and $y$ we obtain the other type of theta characteristic on the curve
\begin{equation*}
(\tilde{\eta},\OO(1))
\end{equation*} 
where $\tilde{\eta}^{\otimes 2}\sim K_{\tilde{C}}$ where the global sections are glued together at the nodes. However, as $h^0(\PP^1,\OO(1))=2$ the values at $x$ and $y$ completely determine the section on $\PP^1$ and we obtain the parity of these theta characteristics is is thus $h^0(\tilde{C},\tilde{\eta})$ mod $2$. There are $2^{g-2}(2^{g-1}+1)$ even and $2^{g-2}(2^{g-1}-1)$ odd theta characteristics of this type. Each have multiplicity $2$ which gives the expected $2^{g-1}(2^g+1)$ even and $2^{g-1}(2^g-1)$ even theta characteristics on $C$.

Consider now a curve $C$ of pseudocompact type with irreducible components $C_1,...,C_k$. Blowing up and inserting an exceptional component at every separating node we obtain the theta characteristic on $C$ to be
\begin{equation*}
(\eta_1,...,\eta_k,\{\OO(1)\}_{i=1}^{k-1})
\end{equation*}
where $\eta_i$ is a theta characteristic on $C_i$ and $\OO(1)$ is a line bundle of degree one on the exceptional $\PP^1$ components. This gives the total degree $\sum_{i=1}^k(g_i-1)+(k-1)=g-1$ as expected and we observe the parity to be \begin{equation*}
\sum_{i=1}^k h^0(C_i,\eta_i) \text{   mod } 2
\end{equation*}
where if any component $C_i$ has self nodes then the $\eta_i$ is of the types discussed earlier.

\subsection{De Jonquieres' formula}\label{sec:dJ}
We will require some tools for enumerating the occurrence of meromorphic differentials of specified signatures in test curves in the moduli space. The first such tool is de Jonquieres' formula which enumerated the number of sections with a specified type of vanishing in a general $g_r^d$. The number of sections with ordered zeros of multiplicity $k_i$ for $i=1,...,\rho$ with $\sum k_i=d$ and $\rho=d-r$ in a general $g_d^r$ on a general genus $g$ curve  is
\begin{equation*}
\dJ[g;k_1,...,k_\rho]=\frac{g!}{(g-\rho-1)!}\prod_{i=1}^{\rho}k_i\left( \sum_{j=0}^{\rho-1}\left(\frac{(-1)^j}{g-\rho+j}\sum_{|I|=j}\left(\prod_{i\nin I}k_i    \right)  \right)+\frac{(-1)^\rho}{g}\right).
\end{equation*}
where $I$ is a subset of $\{1,...,\rho\}$ and $|I|$ denotes the number of elements in $I$. We define 
\begin{equation*}
||I||:=\sum_{i\in I}k_i.
\end{equation*}
This formula can be found throughout the literature in many forms. This presentation is equivalent to that provided in \cite{ACGH} on page $359$ when all $k_i$ are distinct. We present this version for its relative computational ease. It is a variation of that developed in \cite{Coolidge}  page $288$. We will use the convention that $\dJ[1;\emptyset]=1$.

Letting $k_1=r+1$ and $k_i=1$ for $i=2,...,d-r$ recovers the well-known Pl\"ucker formula enumerating the number of simple ramification points in a general $g^r_d$ as
\begin{equation*}
(r+1)d+(r+1)r(g-1)
\end{equation*}
after allowing for the factor $(d-r-1)!$ labelling the simple zeros.

\subsection{The Picard variety method}\label{pic}
De Jonquieres' formula will enumerate the holomorphic sections of general line bundle of a specified type. In some cases, however, we will want to enumerate the of meromorphic sections of a specified type. The Picard variety method enumerates the solutions to a particular equation in the Picard group. Here we provide a summary of this method as presented in~\cite{Mullane}.

For $L$ a specified line bundle of order $d=\sum_{i=1}^gk_i$ on a general genus $g$ curve $C$ we want to enumerate the $(p_1,...,p_g)\in C^g$ such that
\begin{equation*}
\sum_{i=1}^gk_ip_i\sim L.
\end{equation*}
We will follow the treatment in genus $g=2$ of \cite{ChenTarasca} \S2. Consider the map
\begin{eqnarray*}
\begin{array}{cccccc}
f:C^g&\longrightarrow &\Pic^d(C)\\
(p_1,...,p_g)&\longmapsto&\sum_{i=1}^gk_ip_i.
\end{array}
\end{eqnarray*}
The fibre of this map above $L\in\Pic^d(C)$ will give us precisely the solutions of interest. We observe that the domain and range of $f$ are both of dimension $g$. Hence our answer will come from the degree of the map $f$ and an analysis of this fibre. Take a general point $e\in C$ and consider the isomorphism
\begin{eqnarray*}
\begin{array}{cccccc}
h:\Pic^d(C)&\longrightarrow &J(C)\\
L&\longmapsto&L\otimes\OO_C(-de).
\end{array}
\end{eqnarray*}
Now let $F=h\circ f$. Then we have $\deg F=\deg f$. We observe
\begin{equation*}
F(p_1,...,p_g)=\OO_C\biggl(\sum_{i=1}^gk_i(p_i-e)\biggr).
\end{equation*}
Let $\Theta$ be the fundamental class of the theta divisor in $J(C)$. By \cite{ACGH} \S1.5 we have
\begin{equation*}
\deg \Theta^g=g!
\end{equation*}
and the dual of the locus of $\OO_C(k(x-e))$ for varying $x\in C$ has class $k^2\Theta$ in $J(C)$. Hence
\begin{eqnarray*}
\deg F&=&\deg F_*F^*([\OO_C])\\
&=&\deg \left(\prod_{i=1}^g k_i^2\Theta\right)\\
&=&g!\left(\prod_{i=1}^g k_i^2\right)
\end{eqnarray*}
In practice we may want to discount this number by any specific solutions that we may be omitting for some reason. For example, we will be omitting any solutions where $p_i=p_j$ for $i\ne j$. In this case we will need to know not only the existence of any specific solutions that we are discounting by, but also the multiplicity of these solutions. We calculate the multiplicity by investigating the branch locus of $F$. First we look locally analytically at $F$ around each point. If $f_0d\omega,...,f_{g-1}d\omega$ is a basis for $H^0(C,K_C)$, then locally analytically the map becomes
\begin{eqnarray*}
(p_1,...,p_g)&\longmapsto&\biggl(\sum_{i=1}^gk_i\int_e^{p_i}f_0d\omega,...,\sum_{i=1}^gk_i\int_e^{p_i}f_{g-1}d\omega\biggr)
\end{eqnarray*}
modulo $H_1(C,K_C)$. The map on tangent spaces at any fixed point $(p_1,...,p_g)\in C^g$ is the Jacobian of $F$ at the point, which is
\begin{equation*}
\text{DF}(p_1,...,p_g)=\text{diag}(k_1,...,k_g)\begin{pmatrix}f_0(p_1)&...&f_0(p_g)\\
f_1(p_1)&...&f_1(p_g)\\
...&...&...&\\
f_{g-1}(p_1)&...&f_{g-1}(p_g)   \end{pmatrix}
\end{equation*}
Ramification in the map $F$ occurs when the map on tangent spaces is not injective which takes place at the points where $\rk(\text{DF})<g$. The ramification index at a point $(p_1,...,p_g)\in C^g$ will be equal to the vanishing order of the determinant of $\text{DF}{(p_1,...,p_g)}$ at the point. 

We observe that there are two components to the branch locus of $F$. 
\begin{eqnarray*}
\Delta&=& \{(p_1,...,p_g)\in C^g\hspace{0.15cm}|\hspace{0.15cm}p_i=p_j \text{ for some $i\ne j$}\}      \\
\mathcal{K}&=&\{(p_1,...,p_g)\in C^g\hspace{0.15cm}|\hspace{0.15cm}h^0(C,K_C-p_1-...-p_g)>0\}
\end{eqnarray*}
where $\mathcal{K}$ is irreducible and $\Delta$ has $g(g-1)/2$ irreducible components defined by
\begin{equation*}
\Delta_{i,j}= \{(p_1,...,p_g)\in C^g\hspace{0.15cm}|\hspace{0.15cm}p_i=p_j \}
\end{equation*}
for $i,j=1,...,g$ and $i< j$. Hence finding the order of any point in the branch locus will simply be a matter of investigating how these loci meet at the particular point.

\subsection{Divisor theory on $\Mgnb$}
Let $\lambda$ denote the first Chern class of the Hodge bundle on $\Mgnb$ and $\psi_i$ denote the first Chern class of the cotangent bundle on $\Mgnb$ associated to the $i$th marked point where $1\leq i\leq n$. These classes are extensions of classes defined on $\Mgn$ that generates $\Pic(\Mgn)\otimes\QQ$, however, $\Pic(\Mgnb)\otimes\QQ$ contains more classes.

The boundary $\Delta=\Mgnb-\Mgn$ of $\Mgnb$ parameterising marked stable curves of genus $g$ with at least one node is codimension one. Let $\Delta_0$ be the locus of curves in $\Mgnb$ with a non-separating node. Let $\Delta_{i:S}$ for $0\leq i\leq g$, $S\subseteq \{1,...,n\}$ be the locus of curves with a separating node that separates the curve into a genus $i$ component containing the marked points from $S$ and a genus $g-i$ component containing the marked points from $S^C$, the complement of $S$. We require $|S|\geq 2$ for $i=0$ and $|S|\leq n-2$ for $i=g$. We observe that $\Delta_{i:S}=\Delta_{g-i:S^C}$. These boundary divisors are irreducible and can intersect each other and self-intersect. Denote the class of $\Delta_{i:S}$ in $\Pic(\Mgnb)\otimes\QQ$ by $\delta_{i:S}$ and in the case that $n=1$ we denote $\delta_{i:\{1\}}$ by $\delta_i$. See \cite{AC,HarrisMorrison} for more information.

For $g\geq 3$ the these divisor classes freely generate $\Pic(\Mgnb)\otimes \QQ$. For $g=2$ the classes $\lambda,\delta_0$ and $\delta_1$ generate $\Pic(\overline{\mathcal{M}}_2)\otimes \QQ$ with the relation
\begin{equation*}
\lambda=\frac{1}{10}\delta_0+\frac{1}{5}\delta_1.
\end{equation*}
Similarly, $\Pic(\overline{\mathcal{M}}_{2,n})\otimes \QQ$ is freely generated by $\lambda,\psi_i$ and $\delta_{i:S}$ with this relation pulled back under the map $\varphi:\overline{\mathcal{M}}_{2,n}\longrightarrow \overline{\mathcal{M}}_{2}$ that forgets the $n$ marked points.
\begin{equation*}
\lambda=\frac{1}{10}\delta_0+\frac{1}{5}\sum_{1\in S}\delta_{1:S}.
\end{equation*}

\subsection{Maps between moduli spaces}
\label{sec:maps}
There are a number of maps between moduli spaces that prove very useful in divisor class calculations. In this section we present these maps and how the generators of the Picard group pullback under these maps. These results are produced in~\cite{AC} page 161.

Let $(X,q,q_1...,q_{j+1})$ be a general genus $h$ curve marked at $j+2$ general points for $h\geq 0$ and $j\geq -1$ with $j\geq 1$ if $h=0$. Define the map
\begin{eqnarray*}
\begin{array}{cccc}
\pi^{n,j}_{g,h}:&\overline{\mathcal{M}}_{g,n}&\rightarrow& \overline{\mathcal{M}}_{g+h,n+j}\\
&(C,p_1,...,p_n)&\mapsto&(C\bigcup_{p_1=q}X,q_1,p_2,...,p_n,q_2,q_3,...,q_{j+1}).
\end{array}
\end{eqnarray*} 
Letting $\pi=\pi^{n,j}_{g,h}$ for ease of notation, we have for $j\geq 0$
\begin{equation*}
\pi^*\lambda=\lambda,\hspace{0.3cm}\pi^*\delta_0=\delta_0,\hspace{0.3cm}\pi^*\delta_{h:\{1,n+1,...,n+j\}}=-\psi_1.
\end{equation*}
and
\begin{equation*}
\pi^*\psi_i=\begin{cases}  0&\text{  for $i=1$ and $n+1\leq i\leq n+j$}\\ 
\psi_i
  &\text{ otherwise}
\end{cases}
\end{equation*}
Now let $T=\{1,n+1,...,n+j\}$. For $i>h$ 
\begin{equation*}
\pi^*\delta_{i:S}=\begin{cases}  
\delta_{i:S}&\text{  for $S\cap T=\emptyset$}\\
\delta_{i-h:(S\setminus T)\cup\{1\}}&\text{  for $T\subset S$}\\  
0&\text{ otherwise}.\end{cases}
\end{equation*}
For $i=h$ and $h<g$ 
\begin{equation*}
\pi^*\delta_{i:S}=\begin{cases}  
\delta_{i:S}&\text{  for $S\cap T=\emptyset$}\\
\delta_{0:(S\setminus T)\cup\{1\}}&\text{  for $T\subset S$ and $T\ne S$}\\  
-\psi_1&\text{  for $T=S$ }\\  
0&\text{ otherwise}.\end{cases}
\end{equation*}
For $i<h$ 
\begin{equation*}
\pi^*\delta_{i:S}=\begin{cases}  
\delta_{i:S}&\text{  for $S\cap T=\emptyset$ and $i<g$}\\
0&\text{ otherwise}.\end{cases}
\end{equation*}
When $j=-1$ let $(X,q)$ a general genus $h$ curve marked at a general point. The map becomes
\begin{eqnarray*}
\begin{array}{cccc}
\pi^{n,-1}_{g,h}:&\overline{\mathcal{M}}_{g,n}&\rightarrow& \overline{\mathcal{M}}_{g+h,n-1}\\
&(C,p_1,...,p_n)&\mapsto&(C\bigcup_{p_1=q}X,p_2,...,p_n).
\end{array}
\end{eqnarray*} 
Again let $\pi=\pi^{n,-1}_{g,h}$ for ease of notation. Then
\begin{equation*}
\pi^*\lambda=\lambda,\hspace{0.3cm}\pi^*\delta_0=\delta_0,\hspace{0.3cm}\pi^*\psi_i=\psi_{i+1}, \hspace{0.3cm}\pi^*\delta_{h:\emptyset}=-\psi_1+\delta_{h:\emptyset},\hspace{0.3cm}\pi^*\delta_{i:S}=\delta_{i:S}+\delta_{i-h:S\cup\{1\}}\hspace{0.3cm}\text{ for $i:S\ne h:\emptyset$}
\end{equation*}
where $\delta_{i:S}=0$ for $i>g$.

We can also create more complicated maps gluing in multiple tails of different genus with different numbers of marked points to suit our needs. We will describe these maps as needed.

The map $\vartheta$ glues the first and second marked points together
\begin{eqnarray*}
\begin{array}{cccc}
\vartheta:&\overline{\mathcal{M}}_{g,n}&\rightarrow& \overline{\mathcal{M}}_{g+1,n-2}\\
&(C,p_1,...,p_n)&\mapsto&(C/\{p_1\sim p_2\},p_3,...,p_n).
\end{array}
\end{eqnarray*} 
We have 
\begin{equation*}
\vartheta^*\lambda=\lambda,\hspace{0.3cm}\vartheta^*\delta_0=\delta_0+\sum_{p_1\in S,p_2\nin S}\delta_{i:S},\hspace{0.3cm}\vartheta^*\delta_{i:S}=\delta_{i:S}+\delta_{i-1:S\cup\{p_1,p_2\}}\hspace{0.3cm}\vartheta^*\psi_i=\psi_{i+2}.
\end{equation*}

The map $\varphi_j$ forgets the $j$th marked point.
\begin{eqnarray*}
\begin{array}{cccc}
\varphi_j:&\overline{\mathcal{M}}_{g,n}&\rightarrow& \overline{\mathcal{M}}_{g,n-1}\\
&(C,p_1,...,p_n)&\mapsto&(C,p_1,...,p_{j-1},p_{j+1},...,p_n)
\end{array}
\end{eqnarray*} 
We have
\begin{equation*}
\varphi_j^*\lambda=\lambda,\hspace{0.3cm}\varphi_j^*\delta_0=\delta_0,\hspace{0.3cm}\varphi_j^*\delta_{h:S}=\delta_{h:S}+\delta_{h:S\cup\{j\}}\hspace{0.3cm}\varphi_j^*\psi_i=\psi_i-\delta_{0:\{i,j\}}.
\end{equation*}
In the case that $n=1$ let $\varphi:\Mgbm\longrightarrow \Mgb$ be the map forgetting the marked point. In this case for $g$ even we have the one exception that $\varphi^*\delta_{g/2}=\delta_{g/2:\{1\}}$.

\subsection{Divisor Notation}\label{notation}
The divisor notation used in this paper differs based whether the signature $\kappa$ is meromorphic or holomorphic. 

\begin{defn}
For $|\kappa|=g-2+n$ if $\kappa>0$ and $|\kappa|=g-1+n$ if $\kappa<0$ write $\kappa$ in the form
\begin{equation*}
\kappa=(k_1,...,k_n,d_1^{\alpha_1},...,d_m^{\alpha_m})
\end{equation*}
where $d_i\ne d_j$ for $i\ne j$.
Then $D^n_\kappa$ for $n\geq 1$ is the divisor in $\Mgnb$ defined by:
\begin{equation*}
D^n_\kappa:=\frac{1}{\alpha_1!...\alpha_m!}\varphi_*\overline{\PPP}(\kappa)
\end{equation*}
where $\varphi$ forgets the last $g-2$ or $g-1$ marked points for $\kappa>0$ or $\kappa\not>0$ respectively. 
\end{defn}

\section{Effective divisors in $\Mgbm$}\label{Mgbm}
This section provides a simple exposition of the techniques that will be developed in more complicated situations in the next sections. We calculate the classes of divisors in $\Mgbm$ that are the closure of loci of points on smooth curves that form poles or zeros of holomorphic or meromorphic differentials of certain signatures. Unlike in~\cite{Mullane}, where the author used test curves to compute the divisor class, in this paper we will primarily use the method of pulling back divisor classes under maps between different moduli spaces of curves. Through our understanding of the degeneration of meromorphic differentials we are able to explicitly describe the components of the pullback of a divisor coming from the strata of meromorphic differentials. The multiplicity of the components can be computed by an application of the Picard variety method or de Jonquieres' formula. Hence by knowing the class of the components of the pullback of an unknown divisor class or identifying an unknown divisor as a component of the pullback of a known class we obtain the class or many of the coefficients of the class of the unknown divisor.

We record the results of this section as the following theorem.
\begin{thm*}[Theorem~\ref{thmMgbm}]
The divisor $D^1_{g-k,k+1,1^{g-3}}$ for $g\geq3$ and $k=1,...,g-1$ is given by
\begin{eqnarray*}
c_\psi&=&\frac{(k+1)(g-k)((k+1)g^2-(k^2+k+1)g-2)}{2}\\
c_\lambda&=& \frac{ (k+1) (4 - 2 g + 10 k - 2 g k + 11 k^2 + 3 k^3)}{2}\\
 c_0&=&\frac{(k+1)^2-(k+1)^4}{6}\\
 c_i&=&\begin{cases}- \frac{(k+1)(i(g-i+1)(g-i)(k+1)+(g-i-k)((k+1)(g-i)^2-(k^2+k+1)(g-i)-2))}{2}&\text{ for $1\leq i\leq g-k$.}\\
 - \frac{(g - i) (k+1) (-3 g + g^2 + 4 i - g i + 3 k - 4 g k + g^2 k + 
   5 i k - g i k + 3 k^2 - 2 g k^2 + 2 i k^2 + k^3)}{2} &\text{ for $g-k+1\leq i\leq g-1$.}
 \end{cases}
 \end{eqnarray*}
The divisor $D^1_{-h,g+h,1^{g-2}}$ is given by
\begin{eqnarray*}
c_\psi&=&\frac{g(g+h+1)(h-1)(h^2+gh+g+1)}{2}\\
c_\lambda&=&\frac{(1 + g + h) (2 - 3 g^2 + 3 g^3 - 2 h - 4 g h + 9 g^2 h - h^2 + 
   9 g h^2 + 3 h^3)}{2}\\
 c_0&=&\frac{(g+h)^2-(g+h)^4}{6}\\
 c_i&=&\frac{(i-g)}{2} (g^3(2h+i)+g^2(5h^2+2i h+2h+i-1)+g(4h^3+h^2(i+4)+2h(i-1)+i)\\
 &&-2h-i+1+h(h+2)(i+h^2))\hspace{0.3cm}\text{ for $1\leq i\leq g-1$.}
\end{eqnarray*}
\end{thm*}

\begin{rem}
Setting $k=0$ in the first formula we recover $(g-2)W$ for $W=D^1_{g,1^{g-2}}$ the known class of the Weierstrass divisor.
\end{rem}

\subsection{The Weierstrass divisor}\label{sec:w}
Cukierman \cite{Cukierman} calculated the Weierstrass divisor to be  
\begin{equation*}
W=\frac{g(g+1)}{2}\psi-\lambda-\sum_{i=1}^{g-1}\frac{(g-i)(g-i+1)}{2}\delta_i.
\end{equation*}
This is the closure in $\Mgbm$ of Weierstrass points and in our notation we have $W=D^1_{g,1^{g-2}}$.

\subsection{The residual divisor}\label{sec:res}
The first generalisation of the Weierstrass divisor in $\Mgbm$ is the closure of the locus of points that are residual to Weierstrass points. That is
\begin{equation*}
R=\overline{\big\{[C,q] \in \Mgm \hspace{0.1cm}|\hspace{0.1cm} h^0(K_C-gp-q)>0\text{ for some $p\in C$}   \big\}}.
\end{equation*}
We call this the residual divisor and in our notation we have $R=D^1_{1,g,1^{g-3}}$.

In this section we will calculate the class of $R$ through the use of maps between moduli spaces of curves and  previously known classes. For $g\geq 4$ consider the map
\begin{eqnarray*}
\begin{array}{cccc}
\pi:&\overline{\mathcal{M}}_{g-1,1}&\rightarrow& \overline{\mathcal{M}}_{g,1}\\
&(C,y)&\mapsto&(C\bigcup_{x=y}E,q).
\end{array}
\end{eqnarray*} 
that glues a general marked elliptic curve $(E,x,q)$ at $x$ to $y$ at $(C,y)\in\overline{\mathcal{M}}_{g-1,1}$ as described in in Section~\ref{sec:maps}. Consider how $R$ pulls back under this map. If $C$ is a smooth curve there are only two ways to choose a codimension one loci $(C,y)\in\overline{\mathcal{M}}_{g-1,1}$ such that $(C\bigcup_{y=x}E,q)$ occurs with $q$ as a limit of residual points. If $y$ is a Weierstrass point, of a general curve $C$ then $q$ is a limit of residual points and further we see that in this case $q$ is a residual point to $g^2$ points $p$ on $E$ that satisfy
\begin{equation*}
-(g+1)x+q+gp\sim\OO_E
\end{equation*}
Secondly, if $C$ is a genus $g-1$ curve that contains an exceptional Weierstrass point, i.e. a point $p$ such that $h^0(C-gp)>0$ then any point $y$ on such a curve $C$ will make $q$ the limit of residual points. Let $\varphi: \overline{\mathcal{M}}_{g-1,1}\rightarrow \overline{\mathcal{M}}_{g-1}$ be the map that forgets the marked point then the above analysis yields the relation
\begin{equation*}
\pi^*R=g^2W+\varphi^*D
\end{equation*}
where $W$ is closure of the locus of Weierstrass points calculated by Cukierman \cite{Cukierman} and $D$ 
is the closure of the locus of curves containing an exceptional Weierstrass point calculated by Diaz \cite{Diaz}, which agrees with our formula in~\cite{Mullane}. In $ \overline{\mathcal{M}}_{g-1}$ this divisor is
\begin{equation*}
D=\frac{g(g+1)(3g^2-3g+2)}{2}\lambda-\frac{g^2(g-1)(g+1)}{6}\delta_0-\sum_{i=1}^{[(g-1)/2]}\frac{g i(g-i-1)(g+1)^2}{2}\delta_i.
\end{equation*}
\begin{rem}
A simple check shows that a general point in any boundary component $\delta_i$ is not included in this pullback and indeed we have found all components.
\end{rem}
Knowing the classes of $W$ and $D$ by the pullback relations given in Section~\ref{sec:maps} we obtain all coefficients of $R$ except the coefficient of $\psi$. A simple test curve created by allowing the marked point to vary in a fixed general curve provides the coefficient of $\psi$. This well-known test curve\footnote{This test curve is presented in detail in Section~\ref{sec:appendix}} has intersection $2g-2$ times the coefficient of $\psi$. We also know that any general curve has $(g+1)g(g-1)(g-2)$ residual points. Hence we have for $g\geq 4$,
\begin{equation*}
R=\frac{g(g+1)(g-2)}{2}\psi+\frac{g(3g^3-3g+2)}{2}\lambda+\frac{g^2-g^4}{6}\delta_0   +\sum_{i=1}^{g-1}\frac{g(i-g)(g^2i+gi-g+i-1)}{2}\delta_i.
\end{equation*}
In Section~\ref{sec:appendix} we laboriously reproduce this result by the different methods of Porteous' formula and test curves and show that this formula extends to the case $g=3$.

\subsection{Divisors from meromorphic strata}\label{sec:meroMgbm}
We can now pullback the residual divisor under the map 
\begin{eqnarray*}
\begin{array}{cccc}
\pi:&\overline{\mathcal{M}}_{g,1}&\rightarrow& \overline{\mathcal{M}}_{g+h,1}\\
&(C,y)&\mapsto&(C\bigcup_{x=y}X,q).
\end{array}
\end{eqnarray*} 
that glues a general marked genus $h\geq2$ curve $(X,x,q)$ at $x$ to $y$ at $(C,y)\in\overline{\mathcal{M}}_{g,1}$ as introduced in Section~\ref{sec:maps}. The divisor $D^1_{-h,g+h,1^{g-2}}$ in $\Mgbm$ is defined to be
\begin{equation*}
D^1_{-h,g+h,1^{g-2}}=\overline{\big\{[C,q] \in \Mgm \hspace{0.1cm}|\hspace{0.1cm} h^0(K_C+hp-(g+h)q)>0\text{ for some $q\in C$ with $p\ne q$}   \big\}}.
\end{equation*}
By our discussion of twisted canonical divisors we know that if $q$ is the limit of residual points, then for $C$ a general genus $g$ curve, $y$ is either a Weierstrass point, or it satisfies 
\begin{equation*}
h^0(K_C-(g+h)y+hp)\geq 1
\end{equation*}
for some point $p\in C$ with $y\ne p$. Further, there is no codimension one condition on the curve $C$ such that every point $y$ makes $q$ a limit of residual points. Then we know that set theoretically the pullback of $R$ is the union of the divisors $W$ and $D^1_{-h,g+h,1^{g-2}}$  and we are left to find the multiplicities. For the multiplicity of the Weierstrass divisor, we are looking for solutions on a general genus $h$ curve $X$ of the form 
\begin{equation*}
(g+h)p+\sum_{j=1}^{h-1}q_j\sim K_X+(g+2)x-q
\end{equation*}
where $x$ is the node and $q$ is the marked point and we have placed them in general position. The Picard variety method gives $(g+h)^2h$ such solutions and we simply need to discount for the unique solution $p=x$ of order $h-1$, where the order is because $(h-2)x+\sum_{j=1}^{h-1}q_j\sim K_X-q$ so the determinant of the Brill-Noether matrix will vanish with order $h-2$ at this point and hence the point has order $h-1$. We now observe that $D^1_{-h,g+h,1^{g-2}}$ has order one as each point provides a unique solution. Hence we have
\begin{equation*}
\pi^*R=((g+h)^2h-(h-1))W + D^1_{-h,g+h,1^{g-2}}.
\end{equation*} 
\begin{rem}
A simple check shows that a general point in any boundary component $\delta_i$ is not included in this pullback and indeed we have found all components.
\end{rem}

Hence by the pullback formula provided in Section~\ref{sec:maps} we have
\begin{eqnarray*}
\pi^*R&=&\frac{g(g+h)((g+h)^2h+(g+h)h-g-1)}{2}\psi + \frac{(g+h)(3(g+h)^3-3(g+h)+2)}{2}\lambda\\
&&+\frac{(g+h)^2-(g+h)^4}{6}\delta_0+\sum_{i=1}^{g-1}\frac{(g+h)(i-g)((g+h)^2(i+h)+(g+h)(i+h)-g+i-1)}{2}\delta_i.
\end{eqnarray*} 
This gives 
\begin{equation*}
D^1_{-h,g+h,1^{g-2}}=c_\psi\psi+c_\lambda\lambda+c_0\delta_0+\sum_{i=1}^{g-1}c_i\delta_i
\end{equation*}
where
\begin{eqnarray*}
c_\psi&=&\frac{g(g+h+1)(h-1)(h^2+gh+g+1)}{2}\\
c_\lambda&=&\frac{(1 + g + h) (2 - 3 g^2 + 3 g^3 - 2 h - 4 g h + 9 g^2 h - h^2 + 
   9 g h^2 + 3 h^3)}{2}\\
 c_0&=&\frac{(g+h)^2-(g+h)^4}{6}\\
 c_i&=&\frac{(i-g)}{2} (g^3(2h+i)+g^2(5h^2+2i h+2h+i-1)+g(4h^3+h^2(i+4)+2h(i-1)+i)\\
 &&-2h-i+1+h(h+2)(i+h^2))\hspace{0.3cm}\text{ for $1\leq i\leq g-1$.}
\end{eqnarray*}

\begin{rem}
As discussed in Section~\ref{strata:abelian} this divisor is irreducible for $g\geq 3$ and $g=2$ for $h$ odd. When $g=2$ and $h$ is even this divisor has two connected components based on spin structure. $D^{1,\textit{odd}}_{-h,h+2}=5W$ for $4|h$ and $D^{1,\textit{even}}_{-h,h+2}=5W$ otherwise.
\end{rem}

\subsection{The remaining residual divisors}\label{sec:rem}
The classes $D^1_{(g-k,k+1,1^{g-3})}$ for $k=1,...,g-2$ complete the calculation of divisors on $\Mgbm$ coming from strata of differentials with only zeros away from the marked point, all but one of which are simple.
We again consider the map
\begin{eqnarray*}
\begin{array}{cccc}
\pi:&\overline{\mathcal{M}}_{g-i,1}&\rightarrow& \overline{\mathcal{M}}_{g,1}\\
&(C,y)&\mapsto&(C\bigcup_{x=y}X,q).
\end{array}
\end{eqnarray*} 
that glues a general marked genus $i\geq2$ curve $(X,x,q)$ at $x$ to $y$ at $(C,y)\in\overline{\mathcal{M}}_{g,1}$ as introduced in Section~\ref{sec:maps}. If we pullback the divisor $D^1_{g-k,k+1,1^{g-3}}$ for any $i<g-k$ we have
\begin{equation*}
\pi^*D^1_{g-k,k+1,1^{g-3}}=(k+1)^2iW^{g-i}+D^1_{g-i-k,k+1,1^{g-i-3}},
\end{equation*}
as on a general genus $i$ curve $X$ with general marked points $x$ and $q$ we have two ways that $q$ can be the limit of points of the form we require. We are considering solutions of the type
\begin{equation*}
(-g+i-2)x+(k+1)t+(g-k)q+\sum_{j=1}^{i-1}q_j\sim K_X
\end{equation*}
for a general curve $X$ of genus $i$ and fixed general points $x$ and $q$. Hence the $t$ are the ramification points of $|K_X+(g-i+2)x-(g-k)q|$ which is a $g^k_{k+i}$ and hence has
\begin{equation*}
(k+1)(k+i)+(k+1)k(i-1)=(k+1)^2i
\end{equation*}
ramification points (alternatively, this can be found by the Picard variety method). There are no solutions with $t=q$ as this would require $q$ to be the ramification points of $|K_X+(g-i+2)x|$ and not a general point. There are no solutions with $t=x$ provided $i\leq g-k$ as this would contradict the assumption that $x$ and $q$ are general.

The other way we can have a limit of the type we require at $q$ is if we have 
\begin{equation*}
(-g+i+k-2)x+(g-k)q+\sum_{j=1}^{i}q_j\sim K_X
\end{equation*}
which has a unique solution as $h^0(K_X-(i+k-g-2)x-(g-k)q)=1$ for general $x$ and $q$. 

Finally, we again see that a general point of any boundary divisor $\delta_i$ is not included in this pullback and hence there are no extra boundary components.

Hence by the pullback relation
\begin{eqnarray*}
\text{Coefficient of $\lambda$}&=&(k+1)^2i(-1)+\frac{(k+1) (4 - 2 (g-i) + 10 k - 2( g-i) k + 11 k^2 + 3 k^3)}{2}\\
&=& \frac{ (k+1) (4 - 2 g + 10 k - 2 g k + 11 k^2 + 3 k^3)}{2}\\
\text{Coefficient of $\delta_0$}&=&\text{Coefficient of $\delta_0$ in $D^{1}_{2,k+1,1^{k-1}}$}\\
&&  \text{ (by setting $i=g-k-2$ as long as $g-k-2>0$)}\\
&=&\frac{(k+1)^2-(k+1)^4}{6}\\
\text{Coefficient of $\delta_i$}&=&-\text{Coefficient of $\psi$ in $\pi^*D^{1}_{g-k,k+1,1^{g-3}}$}\\
&=& -\frac{(k+1)^2i(g-i+1)(g-i)}{2}- \frac{(k+1)(g-i-k)((k+1)(g-i)^2-(k^2+k+1)(g-i)-2)}{2}\\
&=& -\frac{(k+1)(i(g-i+1)(g-i)(k+1)+(g-i-k)((k+1)(g-i)^2-(k^2+k+1)(g-i)-2))}{2}     \\
&&\text{ (for $i\leq g-k$)}
\end{eqnarray*} 
 In the case that $i> g-k$ we have the different relation
\begin{equation*}
\pi^*D^{1}_{g-k,k+1,1^{g-3}}=((k+1)^2i-(k-g+i))W^{g-i}+D^{1}_{g-i-(k+1),k+1,1^{g-i-2}},
\end{equation*}
because $g-i-(k+1)<0$. We observe that the order of the Weierstrass divisor is because we are considering solutions of the type
\begin{equation*}
(-g+i-2)x+(k+1)t+(g-k)q+\sum_{j=1}^{i-1}q_j\sim K_X
\end{equation*}
for a general curve $X$ of genus $i$ and fixed general points $x$ and $q$. Hence the $t$ are the ramification points of $|K_X+(g-i+2)x-(g-k)q|$ which is a $g^k_{k+i}$ and has
\begin{equation*}
(k+1)(k+i)+(k+1)k(i-1)=(k+1)^2i
\end{equation*}
ramification points (alternatively, this can be computed by the Picard variety method). There are no solutions with $t=q$. There are no solutions with $t=x$ unless $k-1-g+i\geq 0$ in which case we have the unique solution
\begin{equation*}
(k-1-g+i)x+\sum_{j=1}^{i-1}q_j\sim K_X-(g-k)q
\end{equation*}
with order $k-g+i$. The other way that $q$ can be a limit of the type we require is if we have
\begin{equation*}
(-g+i+k-1)x+(g-k)q+\sum_{j=1}^{i-1}q_j\sim K_X
\end{equation*}
which for $i>g-k$ has a unique solution. In this case we have that the other ramification point sits on the other component.

But as $g-i-(k+1)\leq-2$ we calculated the class of $D^{1}_{g-i-(k+1),k+1,1^{g-i-2}}$ earlier in Section~\ref{sec:meroMgbm}. Hence we have
\begin{eqnarray*}
\text{Coefficient of $\lambda$}&=&((k+1)^2i-(k-g+i))(-1)\\
&&+
\frac{1}{2}(1 + (k+1)) (2 - 3 (g-i)^2 + 3 (g-i)^3 - 2 (k+i+1-g)\\
&& - 4 (g-i) (k+i+1-g) + 9 (g-i)^2 (k+i+1-g) - (k+i+1-g)^2 \\
&&+  9 (g-i) (k+i+1-g)^2 + 3 (k+i+1-g)^3)\\
&=& \frac{ (k+1) (4 - 2 g + 10 k - 2 g k + 11 k^2 + 3 k^3)}{2}\\
\text{Coefficient of $\delta_0$}&=&\frac{(k+1)^2-(k+1)^4}{6}\\
\text{Coefficient of $\delta_i$}&=&-\text{Coefficient of $\psi$ in $\pi^*D^{1}_{g-k,k+1,1^{g-3}}$}\\
&=& -\frac{((k+1)^2i-(k-g+i))(g-i+1)(g-i)}{2}\\
&&- \frac{(g-i)(k+2)(i+k-g)((i+k+1-g)^2+(g-i)(i+k+1-g)+(g-i)+1)}{2}\\
&=& -\frac{(g - i) (k+1) (-3 g + g^2 + 4 i - g i + 3 k - 4 g k + g^2 k + 
   5 i k - g i k + 3 k^2 - 2 g k^2 + 2 i k^2 + k^3)}{2}     \\
&&\text{ (for $i> g-k$)}
\end{eqnarray*} 
Putting this together we have for $g\geq 3$ and $k=1,...,g-2$ 
\begin{equation*}
D^1_{(g-k,k+1,1^{g-3})}=c_\psi\psi+c_\lambda\lambda+c_0\delta_0+\sum_{i=1}^{g-1}c_i\delta_i
\end{equation*}
where
\begin{eqnarray*}
c_\psi&=&\frac{(k+1)(g-k)((k+1)g^2-(k^2+k+1)g-2)}{2}\\
c_\lambda&=& \frac{ (k+1) (4 - 2 g + 10 k - 2 g k + 11 k^2 + 3 k^3)}{2}\\
 c_0&=&\frac{(k+1)^2-(k+1)^4}{6}\\
 c_i&=&\begin{cases}- \frac{(k+1)(i(g-i+1)(g-i)(k+1)+(g-i-k)((k+1)(g-i)^2-(k^2+k+1)(g-i)-2))}{2}&\text{ for $1\leq i\leq g-k$.}\\
 - \frac{(g - i) (k+1) (-3 g + g^2 + 4 i - g i + 3 k - 4 g k + g^2 k + 
   5 i k - g i k + 3 k^2 - 2 g k^2 + 2 i k^2 + k^3)}{2} &\text{ for $g-k+1\leq i\leq g-1$.}
 \end{cases}
\end{eqnarray*}

\begin{rem}
Setting $k=g-1$ the coefficients match those computed for the residual divisor $R=D^1_{1^{g-2},g}$ in Section~\ref{sec:res}. Setting $k=0$ we obtain $(g-2)W=(g-2)D^1_{g,1^{g-2}}$.
\end{rem}

\begin{rem}\label{rem:22}
This divisor is defined for $g\geq 3$. As discussed in Section~\ref{strata:abelian}, by Kontsevich and Zorich \cite{KontsevichZorich}, this divisor is irreducible in all cases except $g=3$, $k=1$. In this case we have $D^{1}_{2,2}$ contains two irreducible components distinguished by spin structure. 
\begin{equation*}
D^{1}_{2,2}=\overline{\Theta}_{3,1}+\varphi^*\overline{H}
\end{equation*}
where $\overline{\Theta}_{3,1}$ is the closure of the locus of $(C,p)$ where $p$ lies on a bitangent to a quartic plane curve, $\overline{H}$ is the closure of the locus of hyperelliptic curves in $\overline{\mathcal{M}}_{3}$ and $\varphi:\overline{\mathcal{M}}_{3,1}\longrightarrow \overline{\mathcal{M}}_{3}$ simply forgets the marked point. We know 
\begin{eqnarray*}
\overline{\Theta}_{3,1}&=&14\psi+7\lambda-\delta_0-9\delta_1-5\delta_2    \\
\varphi^*\overline{H}&=&  9\lambda-\delta_0-3\delta_1-3\delta_2. 
\end{eqnarray*}
where the class of $\overline{\Theta}_{3,1}$ was calculated by Farkas~\cite{FarkasBN} and the class of $\overline{H}$ is well-known. 
\end{rem}

\subsection{Comparison with Brill-Noether divisors}\label{sec:compare}
Eisenbud and Harris \cite{EisenbudHarrisIrred} showed that the class of the closure of a pointed Brill-Noether divisor can be expressed as $\mu \mathcal{BN}+ \nu W$, where 
\begin{equation*}
\mathcal{BN}=(g+3)\lambda-\frac{g+1}{6}\delta_0-\sum_{i=1}^{g-1}i(g-i)\delta_i
\end{equation*}
is the pullback from $\Mgb$ of the Brill-Noether divisor, $W$ is the Weierstrass divisor and $\mu$ and $\nu$ are real numbers that are both positive if the divisor is effective and not just a virtual class. We observe that such a divisor will satisfy
\begin{equation*}
\mu=-\frac{6c_0}{g+1}
\end{equation*}
and
\begin{equation*}
\nu=\frac{2c_{\psi}}{g(g+1)} =-\frac{6(g+3)}{g+1}c_0-c_\lambda.
\end{equation*}
Hence as divisors coming from the interior will always have $c_\psi\geq 0$ we have the simple coefficient check
\begin{equation*}
2c_\psi+6(g+3)gc_0+g(g+1)c_\lambda= 0.
\end{equation*}
Any class that violates this cannot be the class of a Brill-Noether divisor. No class calculated in this section satisfies this relation (other than the Weierstrass divisor) and hence do not correspond to the class of a Brill-Noether divisor.

\section{Reproducing known divisor classes in $\Mgnb$}\label{known}
Many divisors in $\Mgnb$ coming from the strata of abelian differential have been calculated under different guises. In this section we will efficiently reproduce these classes, providing an exposition of the method of calculation of classes we will employ in the computation of previously unknown classes in later sections.
  
Let $\underline{d}=(d_1,...,d_n)$ be an $n$-tuple of integers satisfying $\sum d_j=g$ with $d_j\geq 1$. Then 
\begin{equation*}
D_{\underline{d}}:=\{[C;p_1,...,p_n]\in \Mgn \hspace{0.2cm}\big| \hspace{0.2cm}h^0(C,d_1p_1+...+d_np_n)\geq 2\}
\end{equation*}
is a divisor in $\Mgn$. Logan \cite{Logan} showed that the closure of $D_{\underline{d}}$ in $\Mgnb$ has class
\begin{equation}\label{logan}
[\overline{D}_{\underline{d}}]=-\lambda+\sum_{j=1}^n\begin{pmatrix}d_j+1\\2  \end{pmatrix}\psi_j
-0\cdot\delta_0-\sum_{\tiny{\begin{array}{cc}i,S\end{array}}}\begin{pmatrix}|d_S-i|+1\\2  \end{pmatrix}\delta_{i:S}
\end{equation}
in $\Pic(\Mgnb)\otimes \QQ$, where $d_S:=\sum_{j\in S}d_j$. 

By an application of Riemann-Roch and Serre duality we see that for such a $\underline{d}=(d_1,...,d_n)$ we have
\begin{equation*}
h^0(C,d_1p_1+...+d_np_n)=h^0(K_C-d_1p_1...-d_np_n)+1.
\end{equation*}
Hence the class described is on the interior of the moduli space a class coming from abelian differentials and we have $\overline{D}_{\underline{d}}$ is equal to $D^n_{(\underline{d},1^{g-2})}$ on  $\Mgn$. As presented in Section~\ref{strata:abelian}, the divisor $D^n_{\underline{d},1^{g-1}}$ is simply
\begin{equation*}
D^n_{\underline{d},1^{g-2}}=\frac{1}{(g-2)!}\varphi_*\overline{\PPP}(\underline{d},1^{g-2})
\end{equation*}
where $\varphi:\overline{\mathcal{M}}_{g,n+g-2}\longrightarrow\overline{\mathcal{M}}_{g,n}$ forgets the last $g-2$ marked points and $\overline{\PPP}(\underline{d},1^{g-2})$ is the closure of the stratum of canonical divisors with signature $(\underline{d},1^{g-2})$ described in Section~\ref{strata:abelian}.

As the divisors $[\overline{D}_{\underline{d}}]$ and $D^n_{\underline{d},1^{g-1}}$ agree on the interior and we are left to wonder if these two classes agree on the boundary or if they differ by some effective boundary component. However, a general point of any boundary component is not included in the closure as defined by Logan. Similarly, this general point is not the result of the limit of abelian differentials for any signature $(\underline{d},1^{g-2})$. Hence these two classes do in fact agree.

\begin{rem}
By defining these divisors as coming from strata of abelian differentials, Section~\ref{strata:abelian} shows that by the results of Kontsevich and Zorich \cite{KontsevichZorich} the divisors  $\overline{D}_{\underline{d}}$ are irreducible for $g\geq2$.
\end{rem}

This construction can be generalised to allow poles at the marked points. Now let $\underline{d}=(d_1,...,d_n)$ be an $n$-tuple of integers satisfying $\sum d_j=g-1$ with at least one $d_j<0$. Then 
\begin{equation*}
D_{\underline{d}}:=\{[C;p_1,...,p_n]\in \Mgn \hspace{0.2cm}\big| \hspace{0.2cm}h^0(C,d_1p_1+...+d_np_n)\geq 1\}
\end{equation*}
is a divisor in $\Mgn$. M\"uller \cite{Muller} showed using Porteous' formula and test curves that the closure of $D_{\underline{d}}$ in $\Mgnb$ has class\footnote{Note that in this formula $S\ne\{1,...,n\}$. In this case the coefficient is found by $S^C=\emptyset\subset S_+$. The condition on $S$ in the formula is separating the cases where all poles lie on the same component.}
\begin{equation}\label{MGZ}
[\overline{D}_{\underline{d}}]=-\lambda+\sum_{j=1}^n\begin{pmatrix}d_j+1\\2  \end{pmatrix}\psi_j
-0\cdot\delta_0-\sum_{\tiny{\begin{array}{cc}i,S\\S\subset S_+  \end{array}}}\begin{pmatrix}|d_S-i|+1\\2  \end{pmatrix}\delta_{i:S}
-\sum_{\tiny{\begin{array}{cc}i,S\\S\nsubset S_+  \end{array}}}\begin{pmatrix}d_S-i+1\\2  \end{pmatrix}\delta_{i:S}
\end{equation}
in $\Pic(\Mgnb)\otimes \QQ$, where $S_+:=\{j\hspace{0.1cm}|\hspace{0.1cm}d_j>0\}$ and $d_S:=\sum_{j\in S}d_j$. Grushevsky and Zakharov \cite{GruZak} reproduced this result using a different method of a systematic set of test curves.

By an application of Riemann-Roch and Serre duality we see that for such a $\underline{d}=(d_1,...,d_n)$ we have
\begin{equation*}
h^0(d_1p_1+...+d_np_n)=h^0(K_C-d_1p_1...-d_np_n).
\end{equation*}
Hence the class described is on the interior of the moduli space a class coming from meromorphic differentials and we have $\overline{D}_{\underline{d}}$ is equal to $D^n_{(\underline{d},1^{g-1})}$ on the interior $\Mgn$. Here the divisor $D^n_{\underline{d},1^{g-1}}$ is
\begin{equation*}
D^n_{\underline{d},1^{g-1}}=\frac{1}{(g-1)!}\varphi_*\overline{\PPP}(\underline{d},1^{g-1})
\end{equation*}
where $\varphi:\overline{\mathcal{M}}_{g,n+g-1}\longrightarrow\overline{\mathcal{M}}_{g,n}$ forgets the last $g-1$ marked points and $\overline{\PPP}(\underline{d},1^{g-1})$ is the closure of the stratum of canonical divisors with signature $(\underline{d},1^{g-1})$ described in Section~\ref{notation}.

Hence we have found that the divisors $[\overline{D}_{\underline{d}}]$ and $D^n_{(\underline{d},1^{g-1})}$ agree on the interior and we are again left to wonder if these two classes agree on the boundary or if they differ by some effective boundary component. For $g\geq 2$ for any boundary component, a general point is not included in the closure as defined by either M\"uller, Grushevsy and Zakharov's method or by the method of degenerating meromorphic differentials. Hence the two divisors are in fact equal. 

\begin{rem}
By defining these divisors as coming from strata of meromorphic differentials the results of Boissy \cite{Boissy} on the number of connected components imply that  $\overline{D}_{\underline{d}}$ is irreducible for $g\geq2$. In the case that $g=1$ Chen and Coskun \cite{ChenCoskun} showed by the innovative use of a pseudo-automorphism of $\overline{\mathcal{M}}_{1,3}$ that for gcd$(d_1,d_2,d_3)=1$ the irreducible effective divisor $[\overline{D}_{\underline{d}}]$ is extremal. These divisors are not proportional and hence this provides infinitely many extremal effective divisors on $\overline{\mathcal{M}}_{1,3}$ showing that the pseudo-effective cone of $\overline{\mathcal{M}}_{1,n}$ is not polygonal for $n\geq 3$. 
\end{rem}

There has been calculation of classes coming from the components of strata when the strata has more than one irreducible components in one isolated case. Farkas~\cite{FarkasBN} calculated the divisor class of the closure of the locus of points in the support of odd theta characteristics in $\Mgbm$. He denoted this divisor $\overline{\Theta}_{g,1}$, in our notation he calculated 
\begin{equation}\label{far}
D^{1,\textit{odd}}_{2^{g-1}}=2^{g-3}((2^g-1)(\lambda+2\psi)-2^{g-3}\delta_0-\sum_{i=1}^{g-1}(2^i+1)(2^{g-i}-1)\delta_i).
\end{equation}
This is the class of the odd component of $D^1_{2^{g-1}}$ which has $2$ components by Konsevich and Zorich~\cite{KontsevichZorich}. The class of the even component is  
\begin{equation}\label{tex}
D^{1,\textit{even}}_{2^{g-1}}=2^{g-3}((2^g+1)\lambda +0\psi -2^{g-3}\delta_0-\sum_{i=1}^{g-1}(2^i-1)(2^{g-i}-1)\delta_i  ).
\end{equation}
In this case the even theta characteristic gives a cover of $\PP^1$. Hence if $\varphi:\Mgbm\longrightarrow \Mgb$ forgets the marked point we have $D^{1,\textit{even}}_{2^{g-1}}=\varphi^*D^{\textit{even}}_{4,2^{g-3}}$ where the divisor $D^{\textit{even}}_{4,2^{g-3}}$ was originally calculated by Teixidor i Bigas~\cite{Teixidor} as the divisor of curves with a vanishing theta-null.

If $C$ is a general genus $g\geq 3$ curve then any $g-1$ points on $C$ define a hyperplane in the canonical embedding of $C$ that intersects the curve at $g-1$ other points. This give an isomorphism
\begin{equation*}
\gamma:C_{g-1}\longrightarrow C_{g-1}=C^{g-1}/S_{g-1}.
\end{equation*}
Let $\Delta$ be the divisor in $C_{g-1}$ 
\begin{equation*}
\Delta:=\{[p_1,...,p_{g-1}]\in C_{g-1}|\hspace{0.2cm}p_1=p_2\}
\end{equation*} 
and $B$ be the curve
\begin{equation*}
B:=\{[p_1,...,p_{g-1}]\in C_g|\hspace{0.2cm}p_1=p_2\text{ and $p_i$ fixed general for $3\leq i\leq g-1$}\}.
\end{equation*}
Then the numerical classes of $B$ cover $\Delta$ and $B\cdot\Delta=5-3g<0$. Hence as $\Delta$ is irreducible it is extremal. But then $\gamma^*\Delta\cdot \gamma^*B=B\cdot\Delta=5-3g<0$ showing $\gamma^*\Delta$ is also extremal. Globalising this construction and pulling back under the finite morphism $\overline{\mathcal{M}}_{g,g-1}\longrightarrow \overline{\mathcal{M}}_{g,g-1}/S_{g-1}$, Farkas and Verra~\cite{FarkasVerraTheta} computed the class of the closure of the anti-ramification locus, the extremal divisor
\begin{equation}\label{fv}
D^{g-1}_{1^{2g-4},2}=-4(g-7)\lambda+4(g-2)\sum_{i=1}^{g-1}\psi_i-2\delta_0+\sum_{i=0}^g\sum_{s=0}^{i-1}c_{i:s}\delta_{i:S}
\end{equation}
for $s=|S|$ with $s=0,...,i-1$ where
\begin{equation*}
c_{i:s}= -(2g-3)s^2+(4gi+2g-10i+1)s -2gi^2+7i^2-2gi-i-2.
\end{equation*}
Note that $c_{0:s}=c_{g:g-s-1}$ for $s\geq 2$.

\subsection{Logan}\label{ReplicatingLogan}
We now replicate the results of Logan by the use of maps between moduli spaces and degeneration of abelian differentials. We begin by inductively calculating the class of $D^g_{(1^{2g-2})}$  to be
\begin{equation*}
D^g_{1^{2g-2}}=-\lambda+\sum_{i=1}^g\psi_i-  0\delta_0-\sum_{i=0}^{g-1}\frac{(|i-|S||+1)(|i-|S||)}{2}\delta_{i:S}
\end{equation*}
in $\overline{\mathcal{M}}_{g,g}$ where $\delta_{g:S}=\delta_{0:S'}$ when $|S'|\geq 2$ and zero otherwise and the coefficient of $\delta_0$ is zero.

By symmetry we observe that $c_{i:S}=c_{i:S'}$ if $|S|=|S'|$. Hence let $c_{i:n}=c_{i:S}$ for all $|S|=n$. Consider the map $\pi: \Mgbm\longrightarrow\overline{\mathcal{M}}_{g,g}$ that glues in a $\PP^1$ at one of $g+1$ distinct marked points as introduced in Section~\ref{sec:maps}. Pulling back the divisor of interest we obtain the known Weierstrass divisor 
\begin{equation*}
\pi^*D^g_{1^{2g-2}}=D^1_{g,1^{g-2}}=W.
\end{equation*}
This gives 
\begin{equation*}
c_\lambda=-1,\hspace{0.8cm}c_0=0,\hspace{0.8cm}c_{g-i:\emptyset}=c_{i:g}=-\frac{(g-i)(g-i+1)}{2}\hspace{0.8cm}\text{and}\hspace{0.8cm} c_{0:\{1,...,g\}}=c_{0:g}=-\frac{g(g+1)}{2}.
\end{equation*}
Now consider the test curve $B_{i,n}$ created by taking a general genus $i$ curve $X$ marked at $n+1$ general points $x,p_1,...,p_n$ and attaching one of these points to a general genus $g-i$ curve $Y$ marked at $g-n$ points $p_{n+1},...,p_g$. The point of attachment $y$ varies in $Y$. Hence for $n\geq i$
\begin{equation*}
B_{i,n}\cdot D^g_{1^{2g-2}}=(2-2(g-i)-(g-n))c_{i:n}+(g-i)c_{i:n+1}=(n-i)(i^2+gn-gi-in-1)
\end{equation*}
where $c_{i:n}=0$ for $n>g$. The equation in the coefficients holds for all values of $n$. The direct intersection is a result of explicitly enumerating the solutions. In the $X$-aspect, as the marked points are general we have $h^0(K_X+(n-i+2)x-\sum_{j=1}^n p_j)=1$. Hence there is a unique set of $i$ points $\{q_j\}$ (up to labelling) such that
\begin{equation*}
(i-n-2)x+\sum_{j=1}^n p_j+\sum_{j=1}^i q_j\sim K_X
\end{equation*}
In the $Y$-aspect the points $y$ are the ramification points of the system $|K_Y-\sum_{j=n+1}^gp_j|$ and a simple application of the Pl\"ucker formula yields the result. This formula inductively gives the remaining boundary coefficients. 

To obtain the remaining coefficient of $\psi_i$ consider the map $\pi: \overline{\mathcal{M}}_{g,g}\longrightarrow\overline{\mathcal{M}}_{g+1,g+1}$ that glues in a general elliptic tail at one of three marked general points as introduced in Section~\ref{sec:maps}. Under this map 
\begin{equation*}
\pi^*D^{g+1}_{1^{2g}}=D^g_{1^{2g-2}}
\end{equation*}
which gives $c_{\psi_i}=c_{1:2}(D^{g+1}_{1^{2g}})=1$.

We now can replicate Equation~\ref{logan} giving the class of $D^n_{\underline{d},1^{g-2}}$ for $\underline{d}=(d_1,...,d_n)$ an $n$-tuple of integers satisfying $\sum d_j=g$ with $d_j\geq 1$ through the use of maps between moduli spaces.  Consider the map $\pi: \overline{\mathcal{M}}_{g,n}\longrightarrow\overline{\mathcal{M}}_{g,g}$ that glues in a $\PP^1$ at one of $d_j+1$ general marked points at the $j$th marked point as introduced in Section~\ref{sec:maps}. Then we have
\begin{equation*}
\pi^*D^g_{1^{2g-2}}=D^n_{\underline{d},1^{g-2}}.
\end{equation*}
As we know the class on the left we are left to understand how this map pulls back the generators. We observe that $\pi^*\lambda=\lambda$, $\pi^*\delta_0=\delta_0$ and $\pi^*\psi_j=0$ unless the $j$th point was the result of a $d_j=1$ in which case $\pi^*\psi_j=1$. Now if $S$ is the labelled $a_j$ points that were glued in at the $j$th point then $\pi^*\delta_{0:S}=-\psi_j$ and $\pi^*\delta_{i:S}=\delta_{i:\{j\}}$. Similarly if $T\subset \{1,...,n\}$ and $S$ is the set of all $a_j$ points glued in at the $j$th point for all $j\in T$ then $\pi^*\delta_{i:S}=\delta_{i,T}$. Note that this relation also holds for $i=0$ if $|T|\geq 2$. All other classes pullback to give zero. 

Hence we immediately see the class of $D^n_{\underline{d},1^{g-2}}$ to be 
\begin{equation*}
-\lambda+\sum_{j=1}^n\frac{d_j(d_j+1)}{2}\psi_j -0\delta_0-\sum_{i=0,S}^{g-1} \frac{(|i-d_S|+1)(|i-d_S|)}{2}\delta_{i:S}
\end{equation*}
where
\begin{equation*}
d_S=\sum_{j\in S}d_j.
\end{equation*}
\begin{rem}
This is essentially how Logan calculated these divisors as he observed 
\begin{equation*}
\varphi_1{}_*(\delta_{0:\{1,2\}}[\overline{D}_{d_1,...,d_n}])=[\overline{D}_{d_1+d_2,d_3,...,d_n}]
\end{equation*}
and proceeded inductively. One advantage of our method is that it is clear how to generalise this calculation.
\end{rem}

\subsection{M\"uller, Grushevsky and Zakharov}\label{ReplicatingMuller} From our perspective, Equation~\ref{MGZ} generalised the results of Equation~\ref{logan} to allow poles. We replicate this result by using previous results and maps between moduli spaces. Consider now $\underline{d}=(d_1,...,d_m,-h_1,...,-h_{n-m})$ for $n>m$ such that
\begin{equation*}
\sum_{j=1}^md_j   -\sum_{j=1}^{n-m}h_j=g-1.
\end{equation*}
We would like to calculate the class of the divisor $D^n_{\underline{d},1^{g-1}}$ which is the closure in $\Mgnb$ of $[C,q_1,...,q_n]\in\Mgn$ such that 
\begin{equation*}
h^0\bigg(K_C-\sum_{j=1}^md_jq_j   +\sum_{j=1}^{n-m}h_jq_{m+j}\bigg)>0.
\end{equation*}
To generalise our previous method to this situation we just need to understand what to do at the points we require poles. First consider the case that there is exactly one pole, i.e. $m=n-1$. We have $h_1\geq 2$.  Consider the map $\pi: \overline{\mathcal{M}}_{g,n}\longrightarrow\overline{\mathcal{M}}_{g+h_1,n}$ that glues in a general genus $h_1$ curve at one of two marked general points at the $n$th marked point and leaves the first $n-1$ marked points unchanged as introduced in Section~\ref{sec:maps}. Then we have
\begin{equation*}
\pi^*D^{n}_{d_1,...,d_{n-1},1,1^{g+h-2}}=D^n_{d_1,...,d_{n-1},-h_1,1^{g-1}}
\end{equation*}
giving the class of $D^n_{\underline{d},1^{g-1}}$ to be
\begin{equation*}
-\lambda+\sum_{j=1}^m\frac{d_j(d_j+1)}{2}\psi_j+ \frac{h_1(h_1-1)}{2}\psi_n  +  0\delta_0-\sum_{i=0,S\subset\{1,...,n-1\}}^{g}c_{i:S}\delta_{i:S}
\end{equation*}
where 
\begin{equation*}
c_{i:S}=\frac{(|i-d_S|+1)(|i-d_S|)}{2} 
\end{equation*}
for any $S\subset \{1,...,n-1\}$ and $\delta_{g:S}=\delta_{0:S'}$ if $|S'|\geq 2$ where 
\begin{equation*}
d_S=\sum_{j\in S}k_j.
\end{equation*}
Now consider the general case where $\underline{d}=(d_1,...,d_m,-h_1,...,-h_{n-m})$ where $n-m\geq 2$. Let $\sum_{j=1}^{n-m}h_j=h$. Consider the map $\pi: \overline{\mathcal{M}}_{g,m+1}\longrightarrow\overline{\mathcal{M}}_{g,n}$ that glues in a $\PP^1$ at one of $n-m+1$ marked general points at the $(m+1)$th marked point and leaves the first $m$ marked points unchanged as introduced in Section~\ref{sec:maps}. We have
\begin{equation*}
\pi^*D^{n}_{d_1,...,d_m,-h_1,...,-h_{n-m},1^{g-1}}=D^n_{d_1,...,d_m,-h,1^{g-1}}
\end{equation*}
and we immediately obtain the coefficients $c_\lambda=-1$, $c_{\psi_i}=(d_i+1)d_i/2$ for $1\leq i\leq m$, $c_0=0$ and 
\begin{equation*}
c_{i:S}=\frac{(|i-d_S|+1)(|i-d_S|)}{2} 
\end{equation*}
for any $S\subset \{1,...,m\}$. We are left to calculate the coefficients of $\psi_i$ for $m+1\leq i\leq n$ and the coefficients of the boundary classes where the poles don't all lie on the same component.

Consider the map $\pi: \overline{\mathcal{M}}_{g,2}\longrightarrow\overline{\mathcal{M}}_{g,n}$ that glues in a $\PP^1$ at one of $n$ general marked points to the first marked point labelling the new marked points as $\{1,...,n\}-\{j\}$ for some $m+1\leq j\leq n$. We have 
\begin{equation*}
\pi^*D^{n}_{d_1,...,d_m,-h_1,...,-h_{n-m},1^{g-1}}=D^2_{g-1+h_j,-h_j,1^{g-1}}
\end{equation*} 
and we see immediately that $c_{\psi_{j+m}}=(h_j-1)h_j/2$. and $c_{0:S}=(d_S+1)d_S/2$. 

Now for some $S$ such that $|S|=j$ and $S$ or $S'$ are not contained in $\{1,...,m\}$ we consider the map $\pi: \overline{\mathcal{M}}_{g-i,n-j+1}\longrightarrow\overline{\mathcal{M}}_{g,n}$ that glues in a general genus $i$ curve at one of $j+1$ marked general points at the $(n-j+1)$th marked point labelling the $j$ other marked points from $S$ and leaves the first $n-j$ marked points unchanged as introduced in Section~\ref{sec:maps}. We observe that
\begin{equation*}
\pi^*D^{n}_{d_1,...,d_m,-h_1,...,-h_{n-m},1^{g-1}}=D^{n-j+1}_{\underline{d}_{S'},d_S-i,1^{g-i-1}}
\end{equation*}
where $\underline{d}_{S'}$ is the vector made up of the entries of $\underline{d}$ that are indexed by $S'$. We immediately see that 
\begin{equation*}
c_{i:S}=\frac{(d_S-i+1)(d_S-i)}{2}.
\end{equation*}

\subsection{Farkas, Teixidor i Bigas}\label{ReplicatingFarkasTexidor}
Consider the map $\pi:\Mgbm\longrightarrow \overline{\mathcal{M}}_{g+h,1}$ that glues in a general genus $h$ curve at one of two marked general points. We have
\begin{equation*}
\pi^*D^{1,\textit{odd}}_{2^{g+h-1}}=(2^h+1)D^{1,\textit{odd}}_{2^{g-1}} +(2^h-1)D^{1,\textit{even}}_{2^{g-1}}.
\end{equation*}
and
\begin{equation*}
\pi^*D^{1,\textit{even}}_{2^{g+h-1}}=(2^h-1)D^{1,\textit{odd}}_{2^{g-1}} +(2^h+1)D^{1,\textit{even}}_{2^{g-1}}.
\end{equation*}
The class of the components is well known in $g=3$ as the closure of the locus of $(C,p)$ where $p$ lies on a bitangent to a plane curve and the closure of the locus of hyperelliptic curves as discussed in Remark~\ref{rem:22}. This gives the coefficients of $\lambda$ and $\delta_0$.
Further, a simple test curve allowing the marked point to vary in a general curve gives $c_{\psi}(D^{1,\textit{odd}}_{2^{g-1}} )=2^{g-2}(2^g-1)$ and we know the coefficient of $\psi$ in $D^{1,\textit{even}}_{2^{g-1}}$ to be zero.  The equations above then yield the remaining coefficients of all $\delta_i$ for $i>0$ for both components and hence reproduce Equation~\ref{far} and Equation~\ref{tex}.

\subsection{Farkas and Verra}\label{ReplicatingFarkasVerra}
Farkas and Verra~\cite{FarkasVerraTheta} computed the class of the closure of the anti-ramification locus $D^{g-1}_{1^{2g-4},2}$ in $\overline{\mathcal{M}}_{g,g-1}$. Here we replicate this calculation and in later sections we'll use this divisor to compute new divisor classes. Due to the symmetry of the divisor class we'll refer to $c_{i:S}$ by $c_{i:s}$ where $s=|S|$.

First consider the map $\pi:\Mgbm\longrightarrow \overline{\mathcal{M}}_{g,g-1}$ that glues in a $\PP^1$-tail at the marked point at one of $g$ marked general points. We have
\begin{equation*}
\pi^*D^{g-1}_{1^{2g-4},2}=D^1_{g-1,2,1^{g-3}}
\end{equation*}
where the class on the right was computed in Section~\ref{sec:rem}. Hence we obtain the coefficients 
\begin{eqnarray*}
c_\lambda&=& -4(g-7)\\
c_0&=&-2\\
c_{i:0}&=&-(2gi^2-7i^2+2gi+i+2)\\
c_{0:g-1}&=&-(2g+1)(g-1)(g-2)
\end{eqnarray*}
To obtain the final coefficients consider the map $\pi:\overline{\mathcal{M}}_{g,g}\longrightarrow \overline{\mathcal{M}}_{g+1,g}$ that glues in a general elliptic curve at one of two marked general points at the first marked point. We have
\begin{equation*}
\pi^*D^{g}_{1^{2g-2},2}=4D^{g-1}_{1^{2g-2}} +\varphi^*D^{g-1}_{1^{2g-4},2}.
\end{equation*}
where $\varphi:\overline{\mathcal{M}}_{g,g}\longrightarrow \overline{\mathcal{M}}_{g,g-1}$ forgets the first marked point. This relation agrees with our calculated coefficients and for $1\in S$ with $s\leq i$ we obtain the relation
\begin{equation*}
c_{i:S}(D^{g}_{1^{2g-2},2})=-2(i-s+1)(i-s)+c_{i:s}.
\end{equation*}
While for $1\nin S$ for $s\leq i-1$ we have the relation
\begin{equation*}
c_{i:S}(D^{g}_{1^{2g-2},2})=-2(i-s)(i-s-1)+c_{i-1:s-1}.
\end{equation*}
Combining these equations by the symmetry of the boundary coefficients gives the inductive formula
\begin{equation*}
c_{i:s}=c_{i-s:s-1}+4(i-s).
\end{equation*}
The known base case of $c_{i-s:0}$ gives 
\begin{equation*}
c_{i:s}=  -(2g-3)s^2+(4gi+2g-10i+1)s-2gi^2+7i^2-2gi-i-2
\end{equation*}
for $i=1,...,g$ and $s=0,...,i-1$. Note that $c_{0:s}=c_{g:g-s-1}$ for $s\geq 2$. This gives Equation~\ref{fv}.

\section{Coupled partition divisors in $\Mgnb$}\label{Coupled}
When $g\geq 2$ the divisor $D^n_{\underline{d},1^{g-1}}$ is irreducible. We consider the divisors $D^n_{\underline{d},2^{g-1}}$ for $g\geq 2$ with $\underline{d}=(d_1,...,d_n)$ for $\sum_i d_i=0$. When $d_i$ are all even there are two components based on even and odd spin structure. We refer to such a partition $(d_1,...,d_n,2^{g-1})$ of $2g-2$ as a \emph{coupled partition}.

\begin{prop}\label{prop:D112gm2}
The class of the divisor $D^2_{1,1,2^{g-2}}$ in $\overline{\mathcal{M}}_{g,2}$ is
\begin{equation*}
D^2_{1,1,2^{g-2}}=2^{g-3}(2^{g+1}\lambda +2^{g-1}(\psi_1+\psi_2) -2^{g-2}\delta_0    -\sum_{i=0}^{g-1}2^{i+1}(2^{g-i}-1)\delta_{i:\{1,2\}}      -\sum_{i=1}^{g-1}2^{g-1}\delta_{i:\{1\}}         )
\end{equation*}
\end{prop}

\begin{proof}
Consider the map $\pi:\Mgbm\longrightarrow \overline{\mathcal{M}}_{g,2}$ that glues in at the marked point a $\PP^1$-tail at one of three marked distinct points. Under this map we have
\begin{equation*}
\pi^*D^2_{1,1,2^{g-2}}=D^1_{2^{g-1}}=D^{1,\textit{odd}}_{2^{g-1}} +D^{1,\textit{even}}_{2^{g-1}}
\end{equation*}
and hence we obtain the coefficients of $\lambda,\delta_0,\delta_{0:\{1,2\}},\delta_{i:\{1,2\}}$ for $i>0$. For the coefficients of $\psi_i$ consider the test curve $B$ defined by taking a general curve $C$ and marking a general point as the second point. Allow the first point to vary in the curve. We have
\begin{equation*}
B\cdot D^2_{1,1,2^{g-2}}=(2g-1)c_{\psi_1}+c_{\psi_2}+c_{0:\{1,2\}}=\dJ[g;1,2^{g-2}]=2^{g-2}((g-2)2^{g-1}+1)
\end{equation*}
Hence by symmetry
\begin{equation*}
c_{\psi_i}=\frac{1}{2g}(2^{g-2}((g-2)2^{g-1}+1)+2^{g-2}(2^g-1))=2^{2g-4}
\end{equation*}
Finally we need to calculate the coefficients of $\delta_{i:\{1\}}$. Consider the map $\pi:\overline{\mathcal{M}}_{g,2}\longrightarrow \overline{\mathcal{M}}_{g+h,2}$ that glues in at one of the marked points a genus $h$ tail at one of two marked general points. Under this map we have
\begin{equation*}
\pi^*D^2_{1,1,2^{g-2}}=  4^hD^2_{1,1,2^{g-h-2}} .
\end{equation*}
The known $\psi_i$ coefficients then complete our calculation.

\end{proof}

\begin{rem}
This formula is known in $g=2,3$ by Equation~\ref{logan} and Equation~\ref{fv} respectively. 
\end{rem}

We now specialise to the two cases where the signature has exactly one pole which has order two.

\begin{prop}\label{11m2}
Let $\underline{d}=(-2,1,1)$, then
\begin{eqnarray*}
D^{3}_{\underline{d},2^{g-1}}&=&2^{g-3}(2^{g+1}\lambda +2^{g+2}\psi_1+2^{g-1}(\psi_2+\psi_3)-2^{g-2}\delta_0 -\sum_{i=0}^{g-1}2^{i+1}(2^{g-i}-1)\delta_{i\{1,2,3\}}\\
&&-\sum_{i=0}^{g-1}2^{i+1}(2^{g-i}+1)\delta_{i:\{2,3\}}
-\sum_{i=0}^{g-1}2^{g-1}(\delta_{i:\{1,2\}}+\delta_{i:\{1,3\}}).
\end{eqnarray*}
\end{prop}

\begin{proof}
Consider the map $\pi: \overline{\mathcal{M}}_{g,3}\longrightarrow \overline{\mathcal{M}}_{g+1,2}$ that glues in at the first marked point a general elliptic tail at a marked general point. We have
\begin{equation*}
\pi^*D^{2}_{1,1,2^{g-1}}=D^{3}_{-2,1,1,2^{g-1}}+3\varphi_1^*D^2_{1,1,2^{g-2}}
\end{equation*}
where $\varphi_1: \overline{\mathcal{M}}_{g,3}\longrightarrow \overline{\mathcal{M}}_{g,2}$ simply forgets the first marked point. The multiplicity $3$ of the second component represents placing one of the unmarked double zeros at a two-torsion point on the elliptic curve. Now $\pi^*D^2_{1,1,2^{g-1}}$ equals
\begin{eqnarray*}
&=&2^{g-2}(2^{g+2}\lambda +2^{g}(\psi_2+\psi_3)+2^{g+1}\psi_1-2^{g-1}\delta_0    
-\sum_{i=0}^{g-1}2^{i+1}(2^{g+1-i}-1)(\delta_{i:\{2,3\}}+\delta_{i-1:\{1,2,3\}})     \\
&& -\sum_{i=1}^{g-1}2^{g}(\delta_{i:\{2\}}+\delta_{i-1:\{1,2\}})         )\\
&=&2^{g-2}(2^{g+2}\lambda +2^{g}(\psi_2+\psi_3)+2^{g+1}\psi_1 -2^{g-1}\delta_0  
-\sum_{i=0}^{g-1}2^{i+2}(2^{g-i}-1)\delta_{i:\{1,2,3\}}
-\sum_{i=0}^{g-1}2^{i+1}(2^{g+1-i}-1)\delta_{i:\{2,3\}}\\
&& -\sum_{i=1}^{g-1}2^{g}(\delta_{i:\{1,3\}}+\delta_{i:\{1,2\}}) )
\end{eqnarray*}
and $3\varphi_1^*D^2_{1,1,2^{g-2}}$ equals
\begin{eqnarray*}
&=&3\cdot 2^{g-3}(2^{g+1}\lambda +2^{g-1}(\psi_2+\psi_3) -2^{g-2}\delta_0 -2^{g-1}(\delta_{0:\{1,2\}}+\delta_{0:\{1,3\}})   -\sum_{i=0}^{g-1}2^{i+1}(2^{g-i}-1)(\delta_{i:\{2,3\}}+\delta_{i:\{1,2,3\}})   \\
&&   -\sum_{i=1}^{g-1}2^{g-1}(\delta_{i:\{1,3\}}+\delta_{i:\{1,2\}})         ).
\end{eqnarray*}
The Proposition follows.
\end{proof}

\begin{prop}\label{2m2}
Let $\underline{d}=(-2,2)$, then
\begin{equation*}
D^{2}_{\underline{d},2^{g-1}}=2^{g-3}(2^{g+1}\lambda+2^{g+2}\psi_1 +2(2^g+1)\psi_2-2^{g-2}\delta_0 -\sum_{i=0}^{g-1}2^{i+1}(2^{g-i}-1)\delta_{i\{1,2\}}-2^{i+1}(2^{g-i}+1)\sum_{i=1}^{g-1}\delta_{i:\{2\}} ),
\end{equation*}
with
\begin{equation*}
D^{2,\textit{odd}}_{\underline{d},2^{g-1}}=2^{g-3}((2^g-1)\lambda+2(2^g-1)\psi_1+0\psi_2-2^{g-3}\delta_0-\sum_{i=0}^{g-1}(2^i+1)(2^{g-i}-1)\delta_{i:\{1,2\}}-\sum_{i=1}^{g-1}(2^i-1)(2^{g-i}+1)\delta_{i:\{2\}})
\end{equation*}
and
\begin{equation*}
D^{2,\textit{even}}_{\underline{d},2^{g-1}}=2^{g-3}((2^g+1)\lambda+2(2^g+1)\psi_1 +2(2^g+1)\psi_2 -2^{g-3}\delta_0-\sum_{i=0}^{g-1}(2^i-1)(2^{g-i}-1)\delta_{i:\{1,2\}} -\sum_{i=1}^{g-1}(2^i+1)(2^{g-i}+1)\delta_{i:\{2\}} ).
\end{equation*}
\end{prop}

\begin{proof}
For $\underline{d}=(-2,2)$ consider the map $\pi: \overline{\mathcal{M}}_{g,2}\longrightarrow \overline{\mathcal{M}}_{g,3}$ that glues in at the second marked point a $\PP^1$-tail at one of $3$ marked general points and labels these points as the second and third marked points. We have
\begin{equation*}
\pi^*D^{3}_{-2,1,1,2^{g-1}}=D^{2}_{-2,2,2^{g-1}},
\end{equation*}
To distinguish the components consider the map $\pi: \overline{\mathcal{M}}_{g,2}\longrightarrow \overline{\mathcal{M}}_{g+1,1}$ that glues in at the first marked point a general elliptic tail at a marked general point. We have
\begin{equation*}
\pi^*D^{1}_{2^{g}}=D^{2}_{-2,2,2^{g-1}}+3\varphi_1^*D^1_{2^{g-1}}
\end{equation*}
where $\varphi_1: \overline{\mathcal{M}}_{g,3}\longrightarrow \overline{\mathcal{M}}_{g,2}$ simply forgets the first marked point. The multiplicity $3$ of the second component represents placing one of the unmarked double zeros at a two-torsion point to the node on the elliptic curve. On the components this becomes
\begin{equation*}
\pi^*D^{1,\textit{odd}}_{2^{g}}=D^{2,\textit{even}}_{-2,2,2^{g-1}}+3\varphi_1^*D^{1,\textit{odd}}_{2^{g-1}}
\end{equation*}
and
\begin{equation*}
\pi^*D^{1,\textit{even}}_{2^{g}}=D^{2,\textit{odd}}_{-2,2,2^{g-1}}+3\varphi_1^*D^{1,\textit{even}}_{2^{g-1}}.
\end{equation*}
\end{proof}

\begin{rem}
Observe what may at first appear to be the curious consequence that $D^{2,\textit{odd}}_{-2,2,2^{g-1}}=\varphi_2^*D^{2,\textit{odd}}_{2^{g-1}}$. Recall our definition
\begin{equation*}
D^{2,\textit{odd}}_{-2,2,2^{g-1}}:=\overline{\{(C,p_1,p_2)\in\mathcal{M}_{g,2} |-p_1+p_2+s_1+...+s_{g-1}\sim\eta_C \text{ for $\eta_C$ odd and $p_1\ne p_2,s_i$  } \}}.
\end{equation*}
If $h^0(\eta_C)=1$ with $s_1'+...+s_{g-1}'\sim\eta_C$ then $h^0(\eta_C+x)\geq 1$ for any $x$ as $s_1'+...+s_{g-1}'+x$ is a section. However, this section does not satisfy our requirements and hence we require $h^0(\eta_C+x)=2$. Riemann-Roch then gives
\begin{equation*}
h^0(\eta_C-x)=1-g+(g-2)+h^0(\eta_C+x)=1
\end{equation*}
which explains this result.
\end{rem}

\begin{rem}
As a check consider the map $\pi: \overline{\mathcal{M}}_{g-i,2}\longrightarrow \overline{\mathcal{M}}_{g,2}$ that glues in at the first marked point a general genus $i$ tail at one of two marked general points. We obtain
\begin{equation*}
\pi^*D^{2,\textit{odd}}_{\underline{d},2^{g-1}}=2^{i-1}(2^i+1)D^{2,\textit{odd}}_{\underline{d},2^{g-i-1}}+2^{i-1}(2^i-1)D^{2,\textit{even}}_{\underline{d},2^{g-i-1}}
\end{equation*}
and
\begin{equation*}
\pi^*D^{2,\textit{even}}_{\underline{d},2^{g-1}}=2^{i-1}(2^i-1)D^{2,\textit{odd}}_{\underline{d},2^{g-i-1}}+2^{i-1}(2^i+1)D^{2,\textit{even}}_{\underline{d},2^{g-i-1}}
\end{equation*}
\end{rem}

\begin{rem}
When $g=2$ the pinch partition and coupled partition divisors coincide and these results agree with the results in Section~\ref{pp}. 
\end{rem}

At this point we provide a simple example of controlling the residues in a meromorphic differential on  $\PP^1$ that will prove important in our following divisor class calculations.
\begin{ex}\label{P1}
Consider a meromorphic differential on $\PP^1$ with poles of order $j$ and $k$ at $0$ and $\infty$ respectively for $2\leq j\leq k$ and zeros at $1$ and $t$ of order $j+k-m-2$ and $m$ respectively for $j-1\leq m\leq j+k-3$. 
The differential is given locally at $0$ by
\begin{equation*}
c\frac{(z-1)^{j+k-m-2}(z-t)^{m}}{z^j}dz
\end{equation*}  
for some constant $c\ne 0$. The residue at $0$ is given by
\begin{equation*}
c(-1)^{k-1}t^{m-j+1}\sum_{i=0}^{j-1}\begin{pmatrix} j+k-m-2\\ i \end{pmatrix}\begin{pmatrix} m\\ j-i-1 \end{pmatrix}t^i.
\end{equation*}
Hence by investigating the polynomial
\begin{equation*}
\sum_{i=0}^{j-1}\begin{pmatrix} j+k-m-2\\ i \end{pmatrix}\begin{pmatrix} m\\ j-i-1 \end{pmatrix}t^i
\end{equation*}
we obtain the number of meromorphic differentials on $\PP^1$ of signature $\kappa=(-j,-k,m,j+k-m-2)$ with zero residue at the poles.

For example, consider the meromorphic differentials on $\PP^1$  of signature $\mu=(2,h,-h,-4)$. From the above discussion we see that the polynomial becomes
\begin{equation*}
\begin{pmatrix} h\\ 1 \end{pmatrix}t+2\begin{pmatrix} h\\ 2 \end{pmatrix}t^2+\begin{pmatrix} h\\ 3 \end{pmatrix}t^3
\end{equation*}
which has two non-zero solutions when $h\geq 3$ and only one solution when $h=2$. When $h=1$ there are no solutions, indeed the residue at a simple pole is necessarily non-zero. 
\end{ex}

\begin{prop}\label{hmh}
The class of the divisor $D^2_{-h,h,2^{g-1}}$ for $h\geq 3$ in $\overline{\mathcal{M}}_{g,2}$ is\begin{equation*}
D^2_{-h,h,2^{g-1}}=2^{g-2}(2^{g+1}\lambda +2^{g-1}h^2\psi_1+2^{g-1}h^2\psi_2-2^{g-2}\delta_0    -\sum_{i=0}^{g-1}2^{i+1}(2^{g-i}-1)\delta_{i:\{1,2\}}      -\sum_{i=1}^{g-1}2^{g-1}h^2\delta_{i:\{2\}}         ).
\end{equation*}
When $h$ is even this divisor has two components with classes
\begin{eqnarray*}
D^{2,\textit{odd}}_{-h,h,2^{g-1}}&=&2^{g-2}((2^g-1)\lambda+\frac{2^g-1}{4}h^2\psi_1+\frac{2^g-1}{4}h^2\psi_2-2^{g-3}\delta_0-2(2^g-1)\delta_{0:\{1,2\}}-\sum_{i=1}^{g-1}(2^i+1)(2^{g-i}-1)\delta_{i,\{1,2\}}\\
&&-\frac{2^g-1}{4}\sum_{i=1}^{g-1}h^2\delta_{i:\{2\}})
\end{eqnarray*}
and
\begin{eqnarray*}
D^{2,\textit{even}}_{-h,h,2^{g-1}}&=&2^{g-2}((2^g+1)\lambda +\frac{2^g+1}{4}h^2\psi_1+\frac{2^g+1}{4}h^2\psi_2 -2^{g-3}\delta_0-0\delta_{0:\{1,2\}}-\sum_{i=1}^{g-1}(2^i-1)(2^{g-i}-1)\delta_{i,\{1,2\}}\\
&&    -\frac{2^g+1}{4}\sum_{i=1}^{g-1}h^2\delta_{i:\{2\}}).
\end{eqnarray*}
\end{prop}

\begin{rem}
Consider
\begin{equation*}
D_\infty=\lim_{h\to\infty}\frac{1}{h^2}D^2_{-h,h,2^{g-1}}=2^{2g-3}(\psi_1+\psi_2-\sum_{i=1}^{g-1}\delta_{i:\{2\}})
\end{equation*}
with 
\begin{equation*}
D^{\textit{odd}}_\infty= \frac{2^g-1}{2^{g+1}}D_\infty \hspace{1cm}\text{ and }\hspace{1cm}D^{\textit{even}}_\infty=\frac{2^g+1}{2^{g+1}}D_\infty. 
\end{equation*}
Then we obtain
\begin{equation*}
D^{2}_{-h,h,2^{g-1}}=D^{2}_{-2,2,2^{g-1}}+\varphi_1^*D^1_{2^{g-1}}+h^2D_\infty.
\end{equation*}
where $\varphi_1:\overline{\mathcal{M}}_{g,2}\longrightarrow\Mgbm$ forgets the first marked point. This relation also holds in the odd and even spin structure components for even $h$.
\end{rem}

\begin{proof}
Consider the map $\pi:\Mgbm\longrightarrow \overline{\mathcal{M}}_{g,2}$ that glues in at the marked point a $\PP^1$-tail at one of three marked distinct points.  We have
\begin{equation*}
\pi^*D^2_{-h,h,2^{g-1}}=2D^1_{2^{g-1}}=2D^{1,\textit{odd}}_{2^{g-1}} +2D^{1,\textit{even}}_{2^{g-1}},
\end{equation*}
for $h\geq3$  Example~\ref{P1} shows that to obtain a zero residue at the node as required by the global residue condition, there are exactly two points to place the unmarked zero of order $2$ on the $\PP^1$-tail if $h\geq3$ and exactly one point if $h=2$. Hence we again obtain the coefficients of $\lambda,\delta_0,\delta_{0:\{1,2\}},\delta_{i:\{1,2\}}$ for $i>0$. For the coefficients of $\psi_i$ consider the test curve $B_i$ defined by taking a general curve $C$ and allowing the $i$th marked point to vary in the curve while fixing the other marked point at a general point. We have
\begin{equation*}
B_1\cdot D^2_{-h,h,2^{g-1}}=(2g-1)c_{\psi_1}+c_{\psi_2}+c_{0:\{1,2\}}=gh^22^{2g-2}-2^{g-1}(2^g-1)=2^{g-2}(2^g gh^2-2^{g+1}+2)
\end{equation*}
by the Picard variety method where the correction term is for the $2^{g-1}(2^g-1)$ solutions where the points are equal. These solutions violate the global residue condition. Each solution has multiplicity one. Similarly
\begin{equation*}
B_2\cdot D^2_{-h,h,2^{g-1}}=c_{\psi_1}+(2g-1)c_{\psi_2}+c_{0:\{1,2\}}=2^{g-2}(2^g gh^2-2^{g+1}+2)
\end{equation*}
Hence
\begin{equation*}
2gc_{\psi_i}=2^{g-2}(2^g gh^2-2^{g+1}+2)+2^{g-1}(2^g-1),
\end{equation*}
giving 
\begin{equation*}
c_{\psi_i}=2^{2g-3}h^2.
\end{equation*}
Finally, we compute the coefficients of $\delta_{i:\{2\}}$ for $i>0$. Consider the map $\pi: \overline{\mathcal{M}}_{g-i,2}\longrightarrow \overline{\mathcal{M}}_{g,2}$ that glues in at the second marked point a general genus $i$ tail at one of two marked general points. We have
\begin{equation*}
\pi^*D^2_{-h,h,2^{g-1}}=4^iD^2_{-h,h,2^{g-i-1}},
\end{equation*}
and similarly
\begin{equation*}
\pi^*D^2_{h,-h,2^{g-1}}=4^iD^2_{h,-h,2^{g-i-1}},
\end{equation*}
which agrees with all of our calculated coefficients and gives the final unknown coefficients
\begin{equation*}
c_{i:\{2\}}=-4^i2^{2(g-i)-3}h^2=-2^{2g-3}h^2
\end{equation*}
Next we need to identify the components when $h$ is even. As discussed in Section~\ref{strata:abelian}, the divisor has two irreducible components in this case corresponding to odd and even spin structure. We use the same procedure to calculate the class of the components. Consider the map $\pi:\Mgbm\longrightarrow \overline{\mathcal{M}}_{g,2}$ that glues in at the marked point a $\PP^1$-tail at one of three marked distinct points.  By our discussion of meromorphic differentials on $\PP^1$ in Example~\ref{P1} we see that there are two points on the $\PP^1$-tail to place the double zero to make the residue at the node zero and hence satisfy the global residue condition. Further, these will give limits of theta characteristics by Section~\ref{theta} and under this map
\begin{equation*}
\pi^*D^{2,\textit{odd}}_{-h,h,2^{g-1}}=2D^{1,\textit{odd}}_{2^{g-1}},\hspace{1cm}\text{and}\hspace{1cm}\pi^*D^{2,\textit{even}}_{-h,h,2^{g-1}}=2D^{1,\textit{even}}_{2^{g-1}}.
\end{equation*}
Hence we again obtain the coefficients of $\lambda,\delta_0,\delta_{0:\{1,2\}},\delta_{i:\{1,2\}}$ for $i>0$. To obtain the coefficients of $\psi_i$ we need to distinguish which intersections with our test curves $B_1$ and $B_2$ belong to which component. We observe that by the Picard variety method, for any fixed theta characteristic $\eta_C$ on a general curve $C$ there are $gh^2/4$ solutions of the type
\begin{equation*}
\frac{h}{2}p_1+\sum_{j=2}^{g}p_j\sim\eta_C+\frac{h}{2}x
\end{equation*}
and 
$gh^2/4$ solutions of the type
\begin{equation*}
-\frac{h}{2}p_1+\sum_{j=2}^{g}p_j\sim\eta_C-\frac{h}{2}x
\end{equation*}
for any fixed general point $x\in C$. We observe that the solutions we discounted by were all odd theta characteristics and hence
\begin{equation*}
B_1\cdot D^{2,\textit{odd}}_{-h,h,2^{g-1}}=(2g-1)c_{\psi_1}+c_{\psi_2}+c_{0:\{1,2\}}=g(\frac{h}{2})^2 2^{g-1}(2^g-1)-2^{g-1}(2^g-1),
\end{equation*}
and similarly
\begin{equation*}
B_2\cdot D^{2,\textit{odd}}_{-h,h,2^{g-1}}=c_{\psi_1}+(2g-1)c_{\psi_2}+c_{0:\{1,2\}}=g(\frac{h}{2})^2 2^{g-1}(2^g-1)-2^{g-1}(2^g-1).
\end{equation*}
This gives
\begin{equation*}
c_{\psi_1}^{\textit{odd}}=c_{\psi_2}^{\textit{odd}}.
=2^{g-4}(2^g-1)h^2
\end{equation*}
Similarly
\begin{equation*}
B_1\cdot D^{2,\textit{even}}_{-h,h,2^{g-1}}=(2g-1)c_{\psi_1}+c_{\psi_2}+c_{0:\{1,2\}}=2^{g-3}(2^g+1)h^2g
\end{equation*}
and
\begin{equation*}
B_2\cdot D^{2,\textit{even}}_{-h,h,2^{g-1}}=c_{\psi_1}+(2g-1)c_{\psi_2}+c_{0:\{1,2\}}=2^{g-3}(2^g+1)h^2g,
\end{equation*}
giving
\begin{equation*}
c_{\psi_1}^{\textit{even}}=c_{\psi_2}^{\textit{even}}=2^{g-4}(2^g+1)h^2.
\end{equation*}
Finally, consider the map $\pi: \overline{\mathcal{M}}_{g,2}\longrightarrow \overline{\mathcal{M}}_{g+j,2}$ that glues in at the second marked point a general genus $j$ tail at one of two marked general points. We have
\begin{equation*}
\pi^*D^{2,\textit{odd}}_{-h,h,2^{g+j-1}}=2^{j-1}(2^j+1)D^{2,\textit{odd}}_{-h,h,2^{g-1}}+2^{j-1}(2^j-1)D^{2,\textit{even}}_{-h,h,2^{g-1}}
\end{equation*}
and
\begin{equation*}
\pi^*D^{2,\textit{even}}_{-h,h,2^{g+j-1}}=2^{j-1}(2^j-1)D^{2,\textit{odd}}_{-h,h,2^{g-1}}+2^{j-1}(2^j+1)D^{2,\textit{even}}_{-h,h,2^{g-1}}.
\end{equation*}
Similarly
\begin{equation*}
\pi^*D^{2,\textit{odd}}_{h,-h,2^{g+j-1}}=2^{j-1}(2^j+1)D^{2,\textit{odd}}_{h,-h,2^{g-1}}+2^{j-1}(2^j-1)D^{2,\textit{even}}_{h,-h,2^{g-1}}
\end{equation*}
and
\begin{equation*}
\pi^*D^{2,\textit{even}}_{h,-h,2^{g+j-1}}=2^{j-1}(2^j-1)D^{2,\textit{odd}}_{h,-h,2^{g-1}}+2^{j-1}(2^j+1)D^{2,\textit{even}}_{h,-h,2^{g-1}}.
\end{equation*}
This agrees with our calculated coefficients and gives the final unknown coefficients
\begin{equation*}
c_{i:\{2\}}^{\textit{odd}}=-2^{g-4}h^2(2^{g}-1)   \hspace{1cm}\text{and}\hspace{1cm}c_{i:\{2\}}^{\textit{even}}=-2^{g-4}h^2(2^{g}+1).
\end{equation*}
\end{proof}

Next we generalise to $\Mgnb$.
\begin{prop}\label{prop:coupled}
Consider $\underline{d}=(d_1,...,d_{n})$ such that $\sum_{j=1}^{n}d_j=0$ with  $\underline{d}^-\ne \{-2\}$, then
\begin{equation*}
D^{n}_{\underline{d},2^{g-1}}=2^{g-2}(2^{g+1}\lambda +2^{g-1}\sum_{j=1}^nd_j^2\psi_j-2^{g-2}\delta_0   
 -\sum_{\tiny{\begin{array}{cc}|d_S|=0\\|S|\ne n  \end{array}}}\sum_{i=0}^{g}2^{g-i+1}(2^i-1)\delta_{i:S}   
 -2^{g-1}\sum_{|d_S|\geq1}\sum_{i=0}^{g-1}d_S^2\delta_{i:S}        
 ).
\end{equation*}
If all $d_j$ are even then
\begin{equation*}
D^{n,\textit{odd}}_{\underline{d},2^{g-1}}=2^{g-2}((2^g-1)\lambda+\frac{2^g-1}{4}\sum_{j=1}^n d_j^2\psi_j-2^{g-3}\delta_0
-\sum_{\tiny{\begin{array}{cc}|d_S|=0\\|S|\ne n  \end{array}}}\sum_{i=0}^{g}(2^i-1)(2^{g-i}+1)\delta_{i:S}   -\frac{2^g-1}{4}\sum_{|d_S|\geq 2}\sum_{i=0}^{g-1}d_S^2\delta_{i:S}         )
\end{equation*}
and 
\begin{equation*}
D^{n,\textit{even}}_{\underline{d},2^{g-1}}=2^{g-2}((2^g+1)\lambda +\frac{2^g+1}{4}\sum_{j=1}^n d_j^2\psi_j -2^{g-3}\delta_0
-\sum_{|d_S|=0}\sum_{i=0}^{g}(2^i-1)(2^{g-i}-1)\delta_{i:S}     -\frac{2^g+1}{4}\sum_{|d_S|\geq 2}\sum_{i=0}^{g-1}d_S^2\delta_{i:S}         ).
\end{equation*}
\end{prop}

\begin{proof}
Consider the map  $\pi: \overline{\mathcal{M}}_{g,1}\longrightarrow \overline{\mathcal{M}}_{g,n}$ that glues in at the marked point a $\PP^1$-tail at one of $n+1$ marked general points. We have
\begin{equation*}
\pi^*D^{n}_{\underline{d},2^{g-1}}=2D^{1}_{2^{g-1}}
\end{equation*}
where there are two points on the $\PP^1$ tail that will make the residue at the node zero and hence satisfy the global residue condition. If all $d_i$ are even this relation also holds on the odd and even spin structure components. This provides the coefficients for $\lambda,\delta_0$ and $\delta_{i:\{1,...,n\}}$.

Now for any $j\in \{1,...,n\}$ consider the map $\pi: \overline{\mathcal{M}}_{g,2}\longrightarrow \overline{\mathcal{M}}_{g,n}$ that glues in at the second marked point a $\PP^1$-tail at one of $n$ marked general points and labels the remaining points from $\{1,...,n\}\setminus\{j\}$. For $|d_j|\geq 3$ we have
\begin{equation*}
\pi^*D^n_{\underline{d},2^{g-1}}=D^2_{d_j,-d_j,2^{g-1}},
\end{equation*}
when $d_j=2$ we have
\begin{equation*}
\pi^*D^n_{\underline{d},2^{g-1}}=D^2_{2,-2,2^{g-1}}+\varphi_2^*D^1_{2^{g-1}},
\end{equation*}
and if $d_j=-2$ we have
\begin{equation*}
\pi^*D^n_{\underline{d},2^{g-1}}=D^2_{-2,2,2^{g-1}}+\varphi_1^*D^1_{2^{g-1}},
\end{equation*}
where $\varphi_i:\overline{\mathcal{M}}_{g,2}\longrightarrow \overline{\mathcal{M}}_{g,1}$ forgets the $i$th point. When $d_j=1$ we have
\begin{equation*}
\pi^*D^n_{\underline{d},2^{g-1}}=2D^2_{1,1,2^{g-2}},
\end{equation*}
and when $d_j=-1$ we have
\begin{equation*}
\pi^*D^n_{\underline{d},2^{g-1}}=2D^2_{1,1,2^{g-2}}.
\end{equation*}
When all $d_i$ are even these relations also hold on the odd and even spin structure components. These relations agree with the previously calculated coefficients and give us the coefficients for $\psi_j$ and $\delta_{i:\{j\}}$.

Now for any $S\subseteq \{1,...,n\}$ with $2\leq |S|\leq n-2$ consider the map $\pi: \overline{\mathcal{M}}_{g,n-|S|+1}\longrightarrow \overline{\mathcal{M}}_{g,n}$ that glues in at the first marked point a $\PP^1$-tail at one of $|S|+1$ marked general points and labels the remaining points from $S$. For $|d_S|\geq 3$ we have
\begin{equation*}
\pi^*D^n_{\underline{d},2^{g-1}}=D^{n-|S|+1}_{d_S,\underline{d}(S^C),2^{g-1}},
\end{equation*}
in fact as $|S^C|\geq 2$ for $\{d_S,\underline{d}(S^C)\}^-\ne \{-2\}$ this relation holds where we use the convention
\begin{equation*}
D^n_{0,\underline{d},2^{g-1}}=\varphi_1^*D^{n-1}_{\underline{d},2^{g-1}}
\end{equation*}
where $\varphi_i:\overline{\mathcal{M}}_{g,n}\longrightarrow \overline{\mathcal{M}}_{g,n-1}$ forgets the $i$th point. When all $d_i$ are even this relation also holds for the odd and even spin structure components.

The final situation to consider is when $\{d_S,\underline{d}(S^C)\}^-= \{-2\}$. For $d_S=-2$ we have $\underline{d}(S^C)=\{1,1\}$ and hence
\begin{equation*}
\pi^*D^n_{\underline{d},2^{g-1}}=D^{3}_{-2,1,1,2^{g-1}}+\varphi_1^*D^{2}_{1,1,2^{g-2}},
\end{equation*}
while if $d_S\ne -2$, necessarily $\underline{d}(S^C)=\{1,-2\}$ and $d_S=1$, hence
\begin{equation*}
\pi^*D^n_{\underline{d},2^{g-1}}=D^{3}_{1,1,-2,2^{g-1}}+\varphi_3^*D^{2}_{1,1,2^{g-2}}.
\end{equation*}
When all $d_j$ are even these relations hold on the odd and even spin structure components. These relations obtain our remaining coefficients.
\end{proof}

\begin{rem}
As a quick check of the formulas for simple poles consider the divisor $D^4_{-1,-1,1,1,2^{g-1}}$ in $\overline{\mathcal{M}}_{g,4}$. Under the map $\pi:\overline{\mathcal{M}}_{g,3}\longrightarrow \overline{\mathcal{M}}_{g,4}$ that glues in a $\PP^1$-tail at one of three distinct points we have
\begin{equation*}
\pi^*D^4_{-1,-1,1,1,2^{g-1}}=D^3_{-2,1,1,2^{g-1}}+\varphi_1^*D^2_{1,1,2^{g-2}}
\end{equation*}
which agrees with our class calculation.
\end{rem}

\section{Pinch partition divisors in $\Mgnb$}\label{pp}
Farkas and Verra~\cite{FarkasVerraTheta} calculated the class of the divisor $D^{g-1}_{1^{2g-4},2}$. This can be generalised to holomorphic and meromorphic strata with the same signature of "pinched" unmarked points. In the holomorphic case we have:

\begin{prop}\label{prop:pincholo}
Consider $\underline{d}=(d_1,...,d_n)$ with $d_i\geq0$ and $\sum d_i=g-1$ we have for $g\geq 3$
\begin{equation*}
D^n_{(\underline{d},1^{g-3},2)}=-4(g-7)\lambda+\sum_{i=1}^{g-1}(2g(d_i+1)-3d_i-5)d_i\psi_i-2\delta_0+\sum_{i=0}^g\sum_{d_S=0}^{i-1}c_{i:S}\delta_{i:S}
\end{equation*}
where
\begin{equation*}
c_{i:S}=(3-2g)d_S^2+(4gi+2g-10i+1)d_S -2gi^2+7i^2-2gi-i-2.
\end{equation*}
for
\begin{equation*}
d_S:=\sum_{i\in S}d_i
\end{equation*}
and $d_S=0,...,i-1$. Note that $c_{0:S}=c_{g:g-d_S-1}$ for $|S|\geq 2$.
\end{prop}

\begin{proof}
Consider the map $\pi:\Mgnb\longrightarrow \overline{\mathcal{M}}_{g,g-1}$ that glues in a $\PP^1$-tail at each maked point where $d_i\geq 2$. At the $i$th marked point glue a $\PP^1$-tail at one of $d_i+1$ general marked points. Clearly
\begin{equation*}
\pi^*D^{g-1}_{(1^{2g-4},2)}=D^n_{(\underline{d},1^{g-3},2)}.
\end{equation*}
When $n=1$ or $g=3$ we have calculated these classes in other sections.
\end{proof}

To investigate the meromorphic case we begin with $\overline{\mathcal{M}}_{g,2}$. 

\begin{prop}\label{pinch2}For $h\geq 3$ and $g\geq 2$ we have
\begin{eqnarray*}
D^2_{(-h,g+h-2,1^{g-2},2)}&=&(26-4g)\lambda+2h(g h-g-h+2)\psi_1+
2(g+h-2)(g^2+(h-2)g-h)\psi_2\\
&&-2\delta_0+\sum_{i=1}^{g-1}c_{i:\{1\}}\delta_{i:\{1\}}+\sum_{i=1}^gc_{i:\emptyset}\delta_{i:\emptyset}
\end{eqnarray*}
where $\delta_{g:\emptyset}=\delta_{0:\{1,2\}}$ and 
\begin{equation*}
c_{i:\{1\}}=-2((g-3)i^2+(2g h-4h-g+4)i+g h^2-h^2-g h+2h)
\end{equation*}
and
\begin{equation*}
c_{i:\emptyset}=-2((g-3)i^2+(g+1)i+1).
\end{equation*}
When $h=2$ we have
\begin{eqnarray*}
D^2_{(-2,g,1^{g-2},2)}&=&(27-4g)\lambda+4g\psi_1+\frac{g(4g^2-g-9)}{2}\psi_2-2\delta_0+\sum_{i=1}^{g-1}c_{i:\{1\}}\delta_{i:\{1\}}+\sum_{i=1}^gc_{i:\emptyset}\delta_{i:\emptyset}
\end{eqnarray*}
where  
\begin{equation*}
c_{i:\{1\}}=\frac{1}{2}((13-4g)i^2+(17-12g)i-8g)
\end{equation*}
and
\begin{equation*}
c_{i:\emptyset}=\frac{1}{2}((13-4g)i^2-(4g+3)i-4).
\end{equation*}
\end{prop}

\begin{proof}
Consider the map $\pi:\overline{\mathcal{M}}_{g,2}\longrightarrow \overline{\mathcal{M}}_{g+h,2}$ gluing at the first marked point a general genus $h$ curve at one of two marked general points. When $h=1$ we have
\begin{equation*}
\pi^*D^2_{1,g-1,1^{g-2},2}= \varphi_1^*D^1_{g-1,1^{g-3},2}+4D^2_{1,g-1,1^{g-2}},
\end{equation*}
but all three classes here are known and agree with this relation. 

For $h\geq 3$ we have
\begin{equation*}
\pi^*D^2_{1,g+h-2,1^{g+h-3},2}= D^2_{-h,g+h-2,1^{g-2},2} +(4h-2)D^2_{-h+1,g+h-2,1^{g-1}}
\end{equation*}
where the multiplicity $4h-2$ comes from a simple application of the Pl\"ucker formula. The first divisor class is known by Proposition~\ref{prop:pincholo} and the third class by Equation~\ref{logan}. This proves the first equation in the Proposition. 

When $h=2$ we have the relation
\begin{equation*}
\pi^*D^2_{1,g,1^{g-1},2}= D^2_{-2,g,1^{g-2},2}+7\varphi_1^*D^1_{g,1^{g-2}}
\end{equation*}
where $\varphi_1:\overline{\mathcal{M}}_{g,2}\longrightarrow \Mgbm$ forgets the first marked point.
The multiplicity of the second term is due to the Picard variety method and represents the number of solutions to the equation
\begin{equation*}
2p_1+p_2\sim K_C-x+2y
\end{equation*}
for fixed general $x$ and $y$ on a general curve $C$ with genus $g(C)=2$ where $p_1,p_2\ne x,y$. The Picard variety method gives $2^2\cdot 1^2\cdot 2-1$ where the discounted solution is when $p_1=y$ and $p_2=x'$ the conjugate point to $x$ under the hyperelliptic involution. This solution has multiplicity one as $x$ and $y$ were fixed general points. The first divisor class is known by Proposition~\ref{prop:pincholo} and the third class is the Weierstrass divisor given in Section~\ref{sec:w}. This proves the second and final equation in the Proposition. 
\end{proof}

\begin{rem}
For $g\geq 3$ or $g=2$ and $h$ odd these divisors are irreducible. In the case that $g=2$ these divisors correspond to coupled partition divisors. In this case the formula from the two perspectives agree and further for $g=2$ and $h$ even the class of the two irreducible components is given in Proposition~\ref{2m2} and Proposition~\ref{hmh}.
\end{rem}

\begin{rem}
We can perform a quick check on the majority of the calculated coefficients of this Proposition. Consider the map $\pi:\Mgbm\longrightarrow \overline{\mathcal{M}}_{g,2}$ that glues in at the marked point a $\PP^1$-tail at one of three marked distinct points. For $h\geq 3$ we have
\begin{equation*}
\pi^*D^2_{-h,g+h-2,1^{g-2},2}=D^1_{g-1,1^{g-3},2}+2D^1_{g,1^{g-2}}
\end{equation*}
and when $h=2$ this becomes
\begin{equation*}
\pi^*D^2_{-2,g,1^{g-2},2}=D^1_{g-1,1^{g-3},2}+D^1_{g,1^{g-2}}.
\end{equation*}
The change in the multiplicity of the Weierstrass divisor denoted $D^1_{g,1^{g-2}}$ here is due to the fact that in this case there is only one position on the $\PP^1$-tail where the double zero makes the residue at the node vanish. These relations agree with the Proposition.

\end{rem}

This result can be extended to the meromorphic case with exactly one pole. 

\begin{prop}\label{prop:meropinch}
Consider $\underline{d}=(d_1,...,d_n)$ with $\sum d_i=g-2$, $d_j\leq -2$ and $d_i\geq0$ for $i\ne j$, then for $d_j\leq -3$,
\begin{eqnarray*}
D^{n}_{\underline{d},1^{g-2},2}&=&(26-4g)\lambda +\sum_{i=1}^n 2d_i((g-1)d_i+g-2)\psi_i-2\delta_0    +\sum_{i=0}^{g-1}c_{i:S}\delta_{i:S}
\end{eqnarray*}
where for $j\nin S$ and $d_S\leq i-1$,
\begin{equation*}
c_{i:S}=(2-2g) d_S^2+2 ( 2g i+g-4i+1) d_S-2 (g i^2+g i-3i^2 +i+1)  
\end{equation*}
and for $d_S\geq i$,
\begin{equation*}
c_{i:S}= (2-2g) d_S^2  +2 (2g i-g-4i+2) d_S -2 (g i^2-3i^2-g i+4i).
\end{equation*}
For $d_j=-2$
\begin{eqnarray*}
D^{n}_{\underline{d},1^{g-2},2}&=&(27-4g)\lambda +4g\psi_j+\sum_{i\ne j} \frac{(4g(d_i+1)-5d_i-9)d_i}{2}\psi_i-2\delta_0    +\sum_{i=0}^{g-1}c_{i:S}\delta_{i:S}
\end{eqnarray*}
where for $j\nin S$ and $d_S\leq i-1$,
\begin{equation*}
c_{i:S}=\frac{1}{2}((5-4g)d_S^2+(8g i+4g-18i+3)d_S-4g i^2-4gi+13i^2-3i-4)
\end{equation*}
and for $d_S\geq i$,
\begin{equation*}
c_{i:S}=\frac{1}{2}((5-4g)d_S^2+(8g i-4g-18i+9)d_S-4g i^2+4g i+13i^2-17i).
\end{equation*}

\end{prop}

\begin{proof}
Let $d_j=-h$. Consider the map $\pi: \Mgnb\longrightarrow\overline{\mathcal{M}}_{g+h,n}$ that glues in at the $j$th marked point a general genus $h$ curve at one of two marked general points. Then for $h\geq 3$ we have
\begin{equation*}
\pi^*D^n_{\underline{d}',1^{g+h-3},2}=D^n_{\underline{d},1^{g-2},2}+(4h-2)D^n_{\underline{d}'',1^{g-1}}
\end{equation*}
where $\underline{d}'$ and $\underline{d}''$ are the vector $\underline{d}$ with the $j$th entry replaced by a $1$ and $-h+1$ respectively. 

When $h=2$ we obtain
\begin{equation*}
\pi^*D^n_{\underline{d}',1^{g-1},2}=D^n_{\underline{d},1^{g-2},2}+7\varphi_j^*D^{n-1}_{\underline{d}'',1^{g-1}}
\end{equation*}
where $\underline{d}'$ and $\underline{d}''$ are the vector $\underline{d}$ with the $j$th entry replaced by a $1$ and the $j$th entry omitted respectively. The multiplicity $7$ is a result of an application of the Picard variety method as discussed in the proof of Proposition~\ref{pinch2}.
\end{proof}

\begin{rem}
The first example of a pinch partition divisor is $D^{g}_{1^{2g-4},2}$ computed by Farkas and Verra~\cite{FarkasVerraU}. One way to view this divisor class is by considering the flat finite map for a general curve $C$ of genus $g\geq 3$
\begin{eqnarray*}
\begin{array}{cccc}
\varphi:&C^{g-1}&\longrightarrow&C_{g-1}=C^{g-1}/S_{g-1}\\
&(p_1,...,p_{g-1})&\mapsto&[q_1,...,q_{g-1}]
\end{array}
\end{eqnarray*} 
where $S_{g-1}$ acts on $C^{g-1}$ by permuting the points and the map is defined as the $q_i$ satisfying
\begin{equation*}
\sum_{i=1}^{g-1}p_i+\sum_{i=1}^{g-1}q_i\sim K_C.
\end{equation*}
That is, in the canonical embedding of a non-hyperelliptic genus $g$ curve any $g-1$ points $\{p_i\}$ will specify a hyperplane that will intersect the curve at $g-1$ other points we denote $\{q_i\}$. If $\Delta=\{[q_1,...,q_{g-1}]\in C_{g-1}|\hspace{0.2cm}q_1=q_2\}$ and $B$ is the curve in $C_{g-1}$ obtained by setting $q_1=q_2$ a varying point the curve $C$ and fixing the remaining $g-3$ points as general points in $C$. The numerical classes of $B$ cover the irreducible divisor $\Delta$. Hence $B\cdot \Delta=5-3g<0$ implies $\Delta$ is extremal in the pseudo-effective cone by the well-known criteria of a covering curve. But as $\gamma$ is flat and finite we have $\gamma^*B$ provides a covering curve for $\gamma^*\Delta$. Globalising this construction gives the extremal divisor $D^{g}_{1^{2g-4},2}$  in $ \overline{\mathcal{M}}_{g,g-1}$.

From this perspective this construction can be generalised to meromorphic differentials. For any $\underline{d}=(d_1,...,d_g)$ with $\sum d_i=g-2$ and $d_i\in \ZZ\setminus\{0\}$ fix a general genus $g\geq 2$ curve $C$ and the map
\begin{eqnarray*}
\begin{array}{cccc}
\varphi_{\underline{d}}:&C^{g-1}&\dashedrightarrow&C_{g-1}\\
&(p_1,...,p_{g-1})&\mapsto&[q_1,...,q_{g-1}]
\end{array}
\end{eqnarray*} 
where this map is defined as the $q_i\in C$ such that
\begin{equation*}
\sum_{i=1}^{g}d_ip_i+\sum_{i=1}^{g}q_i\sim K_C.
\end{equation*}
The locus of indeterminacy of this map is the codimension two locus where $h^0(K_C-\sum d_ip_i)\geq 2$. Unfortunately, pulling back and completing the curve $B$ does not provide negative intersection with $\gamma^*\Delta$ precisely because of this locus. Interestingly, the image of the locus of indeterminacy under the resolution of $\varphi_{\underline{d}}$ is the points $q_i$ that are colinear in the canonical embedding, that is, $h^0(K_C-q_1-...-q_g)\geq 1$. Globalising this we obtain an extremal divisor as it is contracted by the map from $\widetilde{\mathcal{M}}_{g,g}$ to the universal Picard variety of degree $g$ line bundles~\cite{FarkasVerraU}. 
\end{rem}

\section{Appendix: The residual divisor by alternate methods}\label{sec:appendix}
In this section we reproduce the results of Section~\ref{sec:res} by the different and more labour intensive methods of Porteous' formula and test curves. 

\subsection{Appendix: Locating the limits of Weierstrass and residual points on general nodal curves}
Locating the limits of Weierstrass and residual points on general nodal curves will inform our later analysis.

\subsubsection{A disconnecting node with one component of genus $i$}\label{sec:nodal}
Consider the nodal curve obtained by attaching a general genus $g-i$ curve $Y$ at a non-Weierstrass point $y\in Y$ to a general genus $i$ curve $X$ at a non-Weierstrass point $x\in X$. We would like to locate the limits of Weierstrass and residual points on smooth curves degenerating to out nodal curve. 

Twisted canonical divisors on this nodal curve of the type we are considering have either a point $p_1$ of order $g$ on the $X$ or the $Y$ component or sitting on a $\PP^1$-bridge between $x$ and $y$. Let $p_j$ for $j=2,...,g-1$ be the limits of residual points. If $p_1$ occurs on the $X$ component we have in the $X$-aspect
\begin{equation*}
gp_1 -(g-i+1)x+\sum_{j=2}^{i}p_j   \sim K_X
\end{equation*}
which has $g^2i-i$ solutions by the Picard variety method where we have discounted by the unique solution with $p_1=x$ which has order $i$. In the $Y$-aspect we have
\begin{equation*}
(g-i-1)y+\sum_{j=i+1}^{g-1}p_j\sim K_Y
\end{equation*}
which has a unique solution for a general point $y$. Hence we have $(g^2-1)i$ solutions of this type.

If $p_1$ sits on $Y$ by the same argument we have  $(g^2-1)(g-i)$ solutions.

Further, as $x$ and $y$ are general points on general curves it is not possible to have any $p_j$ on a $\PP^1$-bridge between $x$ and $y$. Any $p_j$ for $j=2,...,g-1$ on such a $\PP^1$-bridge would contradict the curves $X$ and $Y$ or the points $x$ and $y$ being general. If $p_1$ lies on a $\PP^1$-bridge the only possibility is that the bridge contains one zero of multiplicity $g$ and poles at the nodes of order $-(i+1)$ and $-(g-i+1)$. By the cross ratio we can set the poles to $0$ and $\infty$ and the zero to $1$. The resulting differential is given locally at $0$ by
\begin{equation*}
c\frac{(1-z)^{g}}{z^{i+1}}dz
\end{equation*}  
for some constant $c\in \CC^*$. The residues at the nodes then cannot be zero and we have found all
\begin{equation*}
(g^2-1)i+(g^2-1)(g-i)=(g+1)g(g-1)
\end{equation*}
sets of Weierstrass and residual points as expected.

As a cross-check, consider two test curves in $\Mgbm$ constructed from the nodal curve we are considering. Attach a general genus $g-i$ curve $Y$ at a non-Weierstrass point $y\in Y$ to a general genus $i$ curve $X$ at a non-Weierstrass point $x\in X$. Let $B_X$ be the test curve formed by allowing the marked point to vary in the $X$ component and let $B_Y$ be the test curve formed by allowing the marked point to vary in the $Y$ component. Consider the Weierstrass divisor calculated by Cukierman \cite{Cukierman} 
\begin{equation*}
W=\frac{g(g+1)}{2}\psi-\lambda-\sum_{i=1}^{g-1}\frac{(g-i)(g-i+1)}{2}\delta_i.
\end{equation*}
This is the closure in $\Mgbm$ of Weierstrass points and hence intersecting this divisor with our test curves should verify the number of limits of Weierstrass points that we have found on each component of our nodal curve. Indeed
\begin{equation*}
B_X\cdot W=(2i-1)c_\psi-c_i+c_{g-i}=(g+1)(g-1)(g-i)
\end{equation*}
and
\begin{equation*}
B_Y\cdot W=(2(g-i)-1)c_\psi+c_i-c_{g-i}=(g+1)(g-1)i
\end{equation*}
which agree with our calculations.

\subsubsection{A non-disconnecting node}
Consider a general genus $g-1$ curve $X$ and identify two general non-Weierstrass points $x$ and $y$ to form a node. We would like to locate the limits of Weierstrass and residual points on smooth curves degenerating to out nodal curve.  If all $p_j$ occur on $X$ 
\begin{equation*}
gp_1+\sum_{j=2}^{g-1}p_j\sim K_X+x+y
\end{equation*}
with $p_j\ne x,y$. By the Picard variety method we have $g^2(g-1)$ such solutions, but we must discount for any of these solutions where $p_j=x$ or $y$. As $K_X+x$ has a base point at $x$ we see that if any $p_j=y$ then this would cause some other $p_i=x$. If $j\ne1$ and $i\ne 1$ then this causes the curve $X$ to have an exceptional Weierstrass point providing a contradiction with our assumption that $X$ is general. If $i$ or $j=1$ then we have that $x$ or $y$ is a Weierstrass point, contradicting our assumptions.

The last possibility is that the limit of Weierstrass points specialises to the node. In this case we can blow-up and consider this case as $p_1$ sitting on a $\PP^1$-bridge between $x$ and $y$. In this case on the $X$ component of our two component curve we would have
\begin{equation*}
(g-i)x+iy+\sum_{j=2}^{g-1}p_j\sim K_X+x+y
\end{equation*}
 for $i=1,...,  g-1$ which is
\begin{equation*}
(g-i-1)x+(i-1)y+\sum_{j=2}^{g-1}p_j\sim K_X.
\end{equation*}
Hence in the canonical embedding of $X$, fixing such multiplicities at $x$ and $y$ specifies a plane and hence a unique solution. On the $\PP^1$-bridge we have a $g^{g-1}_{2g-2}$ that will adhere to the vanishing of sections in $|K_X+x+y|$ at $x $ and $y$. The limit $p_1$ in this situation will thus be the ramification points of the $g^{g-1}_{g}$ created by imposing the vanishing orders at the nodes in $\PP^1$. The Pl\"ucker formula shows that there are $g$ such points and we obtain $g(g-1)$ solutions of this type. Hence we have found all
\begin{equation*}
g^2(g-1)+g(g-1)=(g+1)g(g-1)
\end{equation*}
sets of Weierstrass and residual points as expected.

\subsection{Porteous' formula}
We calculate the $\lambda$ and $\psi$ coefficients of the residual divisor by realising the locus of interest in $\mathcal{M}_{g,1}$ as the points at which a suitably chosen map between vector bundles drops dimension. The calculation of the class then becomes a well treaded computation in the Chow ring. This method is known as Porteous' formula and we follow the treatment and notation of Faber \cite{Faber}. Let $ \mathcal{C}_g^n$ denote the $n$-fold fibre product of $\Mgm$ over $\Mg$. Consider $\pi_2: \mathcal{C}_g^2\rightarrow \mathcal{C}_g^1$ which forgets the last point and $\pi_1: \mathcal{C}_g^2\rightarrow \mathcal{C}_g^1$ which forgets the first point.

Let $\omega_i$ be the line bundle on $\mathcal{C}_g^n$ obtained by pulling back $\omega$ on $\mathcal{C}_g$ on the projection of the $i$th coordinate and denote its class as $K_i$ in Chow.

Let $E=\EE$ be the Hodge bundle and $F$ be the bundle whose fibres are $H^0(K/(K-gp_1-p_2))$. The bundle $F$ will have rank $g+1$. Then we have 
\begin{equation*}
c(F)=(1 + K_2 - g\Delta_2)(1 + K_1)(1 + 2K_1)\cdot\cdot\cdot(1 + gK_1).
\end{equation*}
We have the natural evaluation map
\begin{equation*}
\varphi:\EE\longrightarrow F
\end{equation*}
where $\EE$ is the Hodge bundle. The locus where this map drops dimension is exactly the points $(C,p_1,p_2)$ where $h^0(C,K_C(-gp_1-p_2))>0$. We calculate the class of this locus via Porteous' Formula. We know 
\begin{equation*}
c(-\EE)=c(\EE^\vee)=1-\lambda_1+\lambda_2...+(-1)^g\lambda_g.
\end{equation*}
Hence by Porteous' Formula we have the class $Y$ of $(C,p_1,p_2)$ where $gp_1+p_2$ is special is the locus where the map has rank $\leq g-1$
\begin{equation*}
\Delta_{g+1-(g-1),g-(g-1)}(c(F)/c(\EE))=\Delta_{2,1}(c(F)\cdot c(\EE^\vee))=c(F)\cdot c(\EE^\vee)\big|_2
\end{equation*}
First we observe that
\begin{eqnarray*}
c_1(F)&=&K_2-g\Delta_2+\frac{g(g+1)}{2}K_1\\
c_2(F)&=&\sum_{i=1}^g iK_1K_2-g(\sum_{i=1}^g i)\Delta_2K_1+\sum_{i=1}^{g-1}(\sum_{j=i+1}^g ij)K_1^2\\
&=& \frac{g(g+1)}{2}K_1K_2-\frac{g^2(g+1)}{2}\Delta_2K_1+\frac{(g-1)g(g+1)(3g+2)}{24}K_1^2,
\end{eqnarray*} 
where
\begin{eqnarray*}
\sum_{i=1}^{g-1}(\sum_{j=i+1}^g ij)&=& \sum_{i=1}^{g-1}i\biggl(\frac{g(g+1)}{2}-\frac{i(i+1)}{2}\biggr)\\
&=&\frac{(g-1)g(g+1)(3g+2)}{24}.
\end{eqnarray*}

Hence we have
\begin{equation*}
[Y]=\lambda_2-\lambda_1(K_2-g\Delta_2+\frac{g(g+1)}{2}K_1)+ \frac{g(g+1)}{2}K_1K_2-\frac{g^2(g+1)}{2}\Delta_2K_1+\frac{(g-1)g(g+1)(3g+2)}{24}K_1^2
\end{equation*}
From Faber for $\pi_d:\mathcal{C}_g^d\rightarrow \mathcal{C}_g^{d-1}$ forgetting the last point we have
\begin{eqnarray*}
\pi_{d*}(MD_{i,d})&=&M \\
\pi_{d*}(MK^k_d)&=& M\cdot\pi_*(\kappa_{k-1})
\end{eqnarray*}
where $M$ is a monomial of classes that are pulled back from $\mathcal{C}_g^{d-1}$. To put a class in a form like this there are a few other relations that are useful
\begin{eqnarray*}
D_{i,j}D_{j,d}&=&D_{i,j}D_{i,d}\hspace{0.5cm}\text{ for $i<j<d$,}\\
D_{i,d}^2&=&-K_iD_{i,d}\hspace{0.5cm}\text{ for $i<d$,}\\
K_dD_{i,d}&=&K_iD_{i,d}.
\end{eqnarray*}
Now via Faber's algorithm for pushing down, forgetting the second point we have
\begin{eqnarray*}
W=\frac{1}{g-2}\pi_{2*}[Y]&=&\frac{1}{g-2}((g-\kappa_0)\lambda_1+\frac{g(g+1)}{2}(\kappa_0-g)K_1)\\
&=&\frac{g(g+1)}{2}\psi -\lambda
\end{eqnarray*}
as expected. Here we have used the fact that on $W$ we have $\pi_{2*}$ is of degree $g-2$ (there are $g-2$ residual points for each Weierstrass point $p$). Forgetting the first point we have
\begin{eqnarray*}
R=\pi_{1*}[Y]&=&(g-\frac{g(g+1)}{2}\kappa_0)\lambda_1+\frac{g(g+1)}{2}(\kappa_0-g)K_2+  \frac{(g-1)g(g+1)(3g+2)}{24}\kappa_1     \\
&=&\frac{g(g+1)(g-2)}{2}\psi+\frac{g(3g^3-3g+2)}{2}\lambda.
\end{eqnarray*}
which agrees with our previous calculation of these coefficients in Section~\ref{sec:res}.

\subsection{Test curves}
By creating a number of curves in $\Mgbm$ we can calculate the intersections with the generators of $\Pic(\Mgbm)\otimes\QQ$ and the residual divisor directly to obtain a number of relationships between the coefficients of the residual divisor. With enough relationships we can determine all coefficients. To this end let the class of the residual divisor in $\Pic(\Mgbm)\otimes\QQ$ be denoted
\begin{equation*}
R=c_\psi \psi+c_\lambda \lambda+\sum_{i=0}^{g-1}c_i\delta_i.
\end{equation*}

\subsubsection*{Test curve $A$}
Consider a general genus $g$ curve $C$. Allow the marked point $q$ to vary in the curve. We observe that the intersection of this test curve with all boundary divisors and $\lambda$ is zero and we have
\begin{equation*}
A\cdot R=-(2-2g)c_\psi.
\end{equation*}
To find this intersection directly we observe that as the curve is general there are $(g+1)g(g-1)$ normal Weierstrass points each with $g-2$ residual points. Hence
\begin{equation*}
A\cdot R=(g+1)g(g-1)(g-2)
\end{equation*}
and hence
\begin{equation*}
c_\psi=\frac{(g+1)g(g-2)}{2}.
\end{equation*}
Hence we have verified this result from Porteous' formula via test curves.

\subsubsection*{Test curves $B_i$}
Let $Y$ be a genus $g-i$ curve and $X$ be a genus $i$ curve. Attach $X$ to $Y$ at a general point in $x\in X$ and allow the attaching point $y\in Y$ to vary. Mark a general point $q\in X$ as shown in Figure $4$. We require $0<i<g-1$.

\begin{figure}[htbp]
\begin{center}
\begin{overpic}[width=0.5\textwidth]{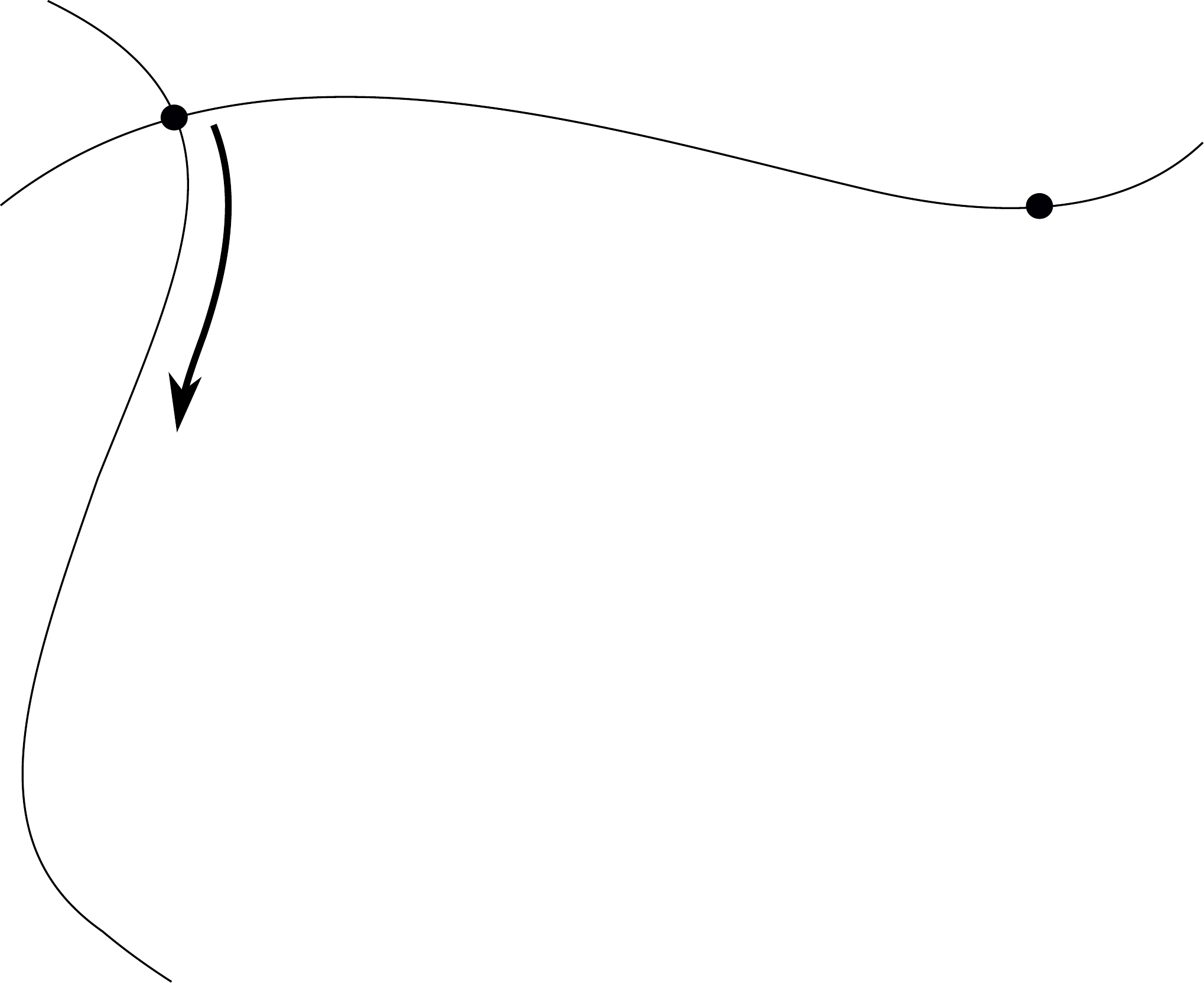}
\put (15, 74){$x$}
\put (65,72){$X$}
\put (85,60){$q$}
\put (45,63){$g(X)=i$}
\put (12,67){$y$}
\put (11,30){$g(Y)=g-i$}
\put (11,18){$Y$}
\put (20,-5){\footnotesize{Figure 4: Test curve $B_i$ }}
\end{overpic}
\end{center}
\end{figure}

We observe that $B_i\cdot R=(2-2(g-i))c_i$. To locate the limits of residual points in this test curve there are two possible cases based on the order of vanishing at $y$. In the first case we have $y$ is a Weierstrass point and in the $Y$-aspect we have
\begin{equation*}
(g-i)y+\sum_{j=1}^{g-i-2}q_j\sim K_Y
\end{equation*}
which has solutions where $y$ is a Weierstrass point and the $q_j$ residual in the $Y$-aspect. There are $(g-i+1)(g-i)(g-i-1)$ such solutions allowing for the ordering of the $q_j$. In the $X$-aspect we have
\begin{equation*}
gp+\sum_{j=g-i-1}^{g-3}q_j\sim K_X+(g-i+2)x-q.
\end{equation*} 
By the Picard variety method we have $g^2i-(i-1)$ solutions where we have discounted for the unique solution where $p=x$ which has order $(i-1)$. Hence solutions of this type contribute 
\begin{equation*}
(g-i+1)(g-i)(g-i-1)(g^2i-(i-1))
\end{equation*}
to the intersection with the residual divisor.

If $y$ is not a Weierstrass point then for any solution we must have in the $Y$-aspect
\begin{equation*}
gp-iy+\sum_{j=1}^{g-i-2}q_j\sim K_Y
\end{equation*}
by the Picard variety method we have 
\begin{equation*}
g^2i^2(g-i)(g-i-1)-(g-i+1)(g-i+1)(g-i)(g-i-1)
\end{equation*}
solutions where we have discounted for the order $g-i+1$ solutions where $p=y$ is a Weierstrass point. We observe that this is consistent with the case $i=1$ where there are no solutions. In the $X$ aspect this corresponds to the unique solution  
\begin{equation*}
\sum_{j=g-i-1}^{g-3}q_j\sim K_X-q-(i-2)x.
\end{equation*}
Hence solutions of this type contribute
\begin{equation*}
(g^2i^2-(g-i+1)^2)(g-i)(g-i-1)
\end{equation*}
to the intersection with the residual divisor and we are left with
\begin{equation*}
B_i\cdot R=g(g^2i+gi-g+i-1)(g-i)(g-i-1)
\end{equation*}
and hence
\begin{equation*}
c_i=-\frac{g(g^2i+gi-g+i-1)(g-i)}{2}
\end{equation*}
for $1\leq i\leq g-2$.

\subsubsection*{Test curves $C_i$}
Let $Y$ be a genus $g-i$ curve and $X$ be a genus $i$ curve. Attach $X$ to $Y$ at a general point in $x\in X$ and a general point $y\in Y$ to vary. Let the marked point $q$ vary in $X$ as shown in Figure $5$. 
\begin{figure}[htbp]
\begin{center}
\begin{overpic}[width=0.5\textwidth]{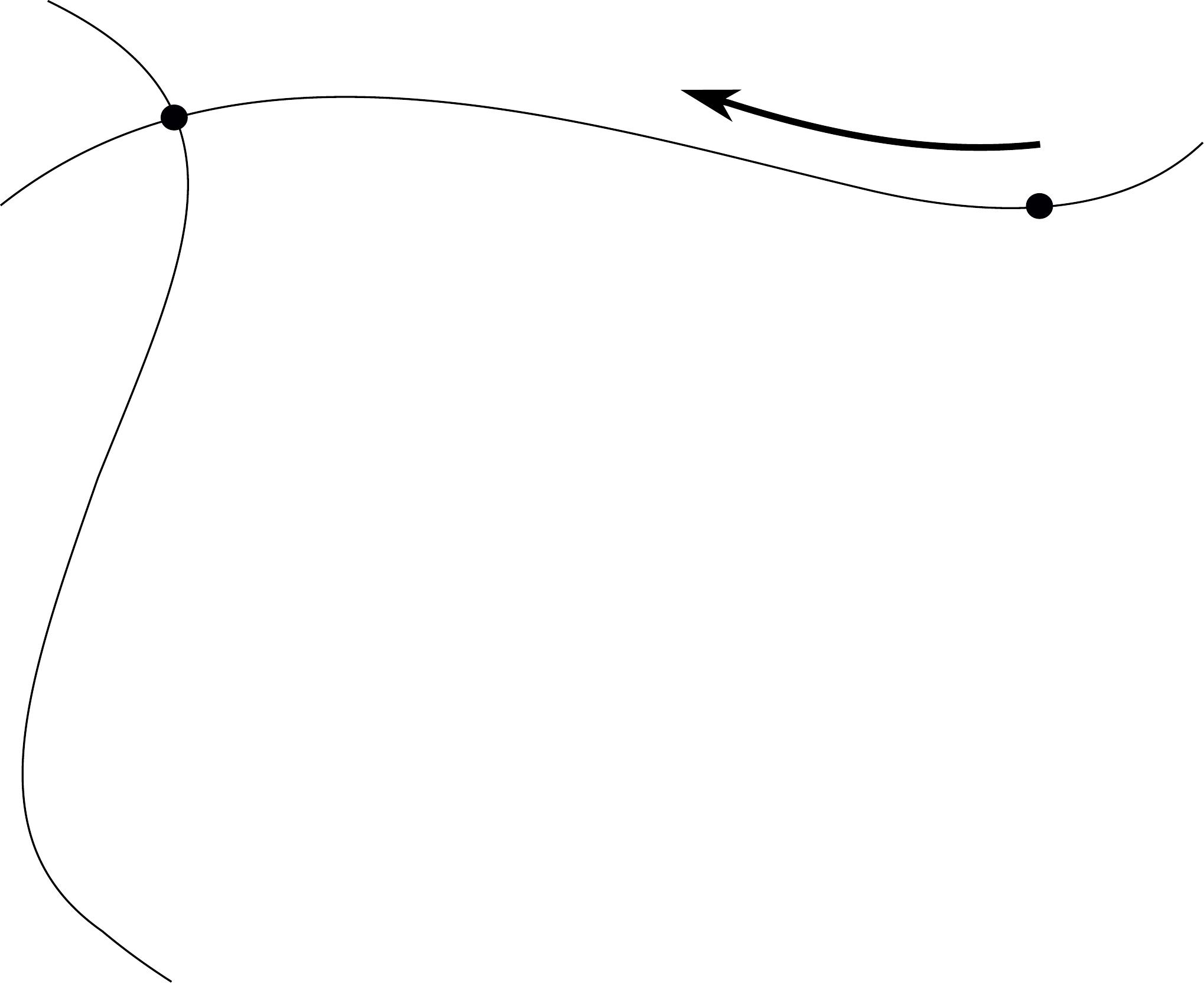}
\put (15, 74){$x$}
\put (35,76){$X$}
\put (85,60){$q$}
\put (45,63){$g(X)=i$}
\put (12,67){$y$}
\put (11,30){$g(Y)=g-i$}
\put (11,18){$Y$}
\put (20,-5){\footnotesize{Figure 5: Test curve $C_i$ }}
\end{overpic}
\end{center}
\end{figure}
We observe that 
\begin{equation*}
C_i\cdot R=(2i-1)c_\psi-c_i+c_{g-i}.
\end{equation*}
But in Section~\ref{sec:nodal} we located the limits of residual points on a general nodal curve of the type we are considering here giving 
\begin{equation*}
C_i\cdot R=g(g^2-1)(i-1)
\end{equation*}
which agrees with our formula for $c_i$ in the last section and shows that it also applies to $i=g-1$.

\subsubsection*{Test curve $D$}
Take a smooth general genus $g-1$ curve $C$. Create a node by identifying one non-special fixed point $y$ on the curve with another point $x$ that varies in the curve. Mark a general point $q$ on the curve as shown in Figure $6$.

\begin{figure}[htbp]
\begin{center}
\begin{overpic}[width=0.25\textwidth]{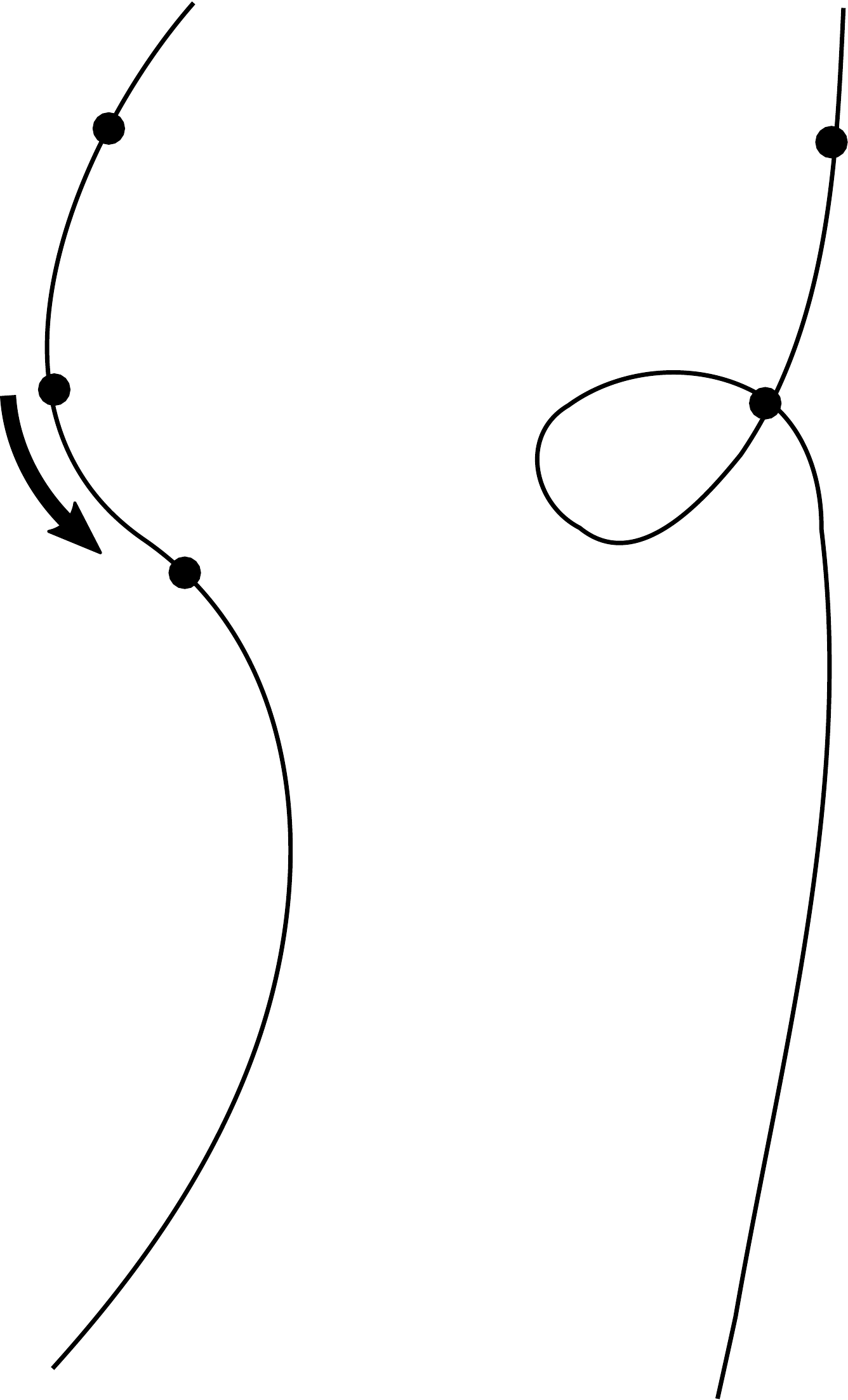}

\put (1, 89){$q$}
\put (53, 89){$q$}
\put (8, 71){$x$}
\put (18,60){$y$}
\put (-15,35){$g(C)=g-1$}
\put (7,18){$C$}
\put (60,68){$x\sim y$}
\put (20,-5){\footnotesize{Figure 6: Test curve $D$ }}
\end{overpic}
\end{center}
\end{figure}
We have
\begin{equation*}
D\cdot R=c_\psi+(2-2g)c_0+c_{g-1}
\end{equation*}
To find the limits of residual points in this test curve there are solutions of two types. The solutions of the first type are of the form
\begin{equation*}
gp+q+\sum_{j=1}^{g-3}q_j\sim K_C+x+y
\end{equation*}
where $y$ and $q$ are fixed and $x$ is varying. The Picard variety method there are $g^2(g-1)(g-2)$ solutions. There are no solutions to discount for with $x=p,q,y,q_j$ or $y=p,q_j$. If $x=p$ then we have
\begin{equation*}
(g-1)x+q+\sum_{j=1}^{g-3}q_j\sim K_C+y.
\end{equation*}
But as $y$ is a base point of $K_C+y$ we have either $y=q_j$ for some $j$ making $q$ a residual point or $y=x$ giving
\begin{equation*}
(g-2)y+q+\sum_{j=1}^{g-3}\sim K_C.
\end{equation*}
Both cases contradict the assumption that $q$ and $y$ are general. Similarly, if $x=q,y,q_j$ for some $j$ or $y=p,q_j$ we have a contradiction of $y$ and $q$ being general points or $C$ a general curve with only normal Weierstrass points.

The second way that limits of residual points can occur in our test curve is if the point approaches the node and $p$ actually sits on a $\PP^1$ between $x$ and $y$. If this occurs we have
\begin{equation*}
ix+(g-i)y+q+\sum_{j=1}^{g-3}q_j\sim K_C+x+y
\end{equation*}
which becomes
\begin{equation*}
(i-1)x+\sum_{j=1}^{g-3}q_j\sim K_C-q-(g-i-1)y
\end{equation*}
for $i=2,...,g-1$. Such $x$ are the ramification points of $|K_C-q-(g-i-1)y|$ which is a $g^{i-2}_{g+i-4}$ and hence by the Pl\"ucker formula there are
\begin{equation*}
(r+1)d+(r+1)r(g(C)-1)=(i-1)(gi-g-i)
\end{equation*}
such ramification points. Hence in total we have
\begin{equation*}
\sum_{i=2}^{g-1}(i-1)(gi-g-i)=\frac{1}{6}g(2g^3-11g^2+19g-10)
\end{equation*}
solutions each with order $g$. 
\begin{rem}
I need to insert an explanation of this multiplicity.
\end{rem}

This gives the relation
\begin{equation*}
c_\psi+(2-2g)c_0+c_{g-1}=g^2(g-1)(g-2)+\frac{1}{6}g^2(2g^3-11g^2+19g-10)
\end{equation*}
and from the known values of $c_\psi$ and $c_{g-1}$ this gives
\begin{equation*}
c_0=\frac{g^2-g^4}{6}.
\end{equation*}

\subsubsection*{Test curve $E$}
Take a pencil of plane cubics. Attach one base point to a general genus $g-1$ curve $C$ at a general point $y\in C$. Mark another general point $q$ on $C$ as shown in Figure $7$.

\begin{figure}[htbp]
\begin{center}
\begin{overpic}[width=0.35\textwidth]{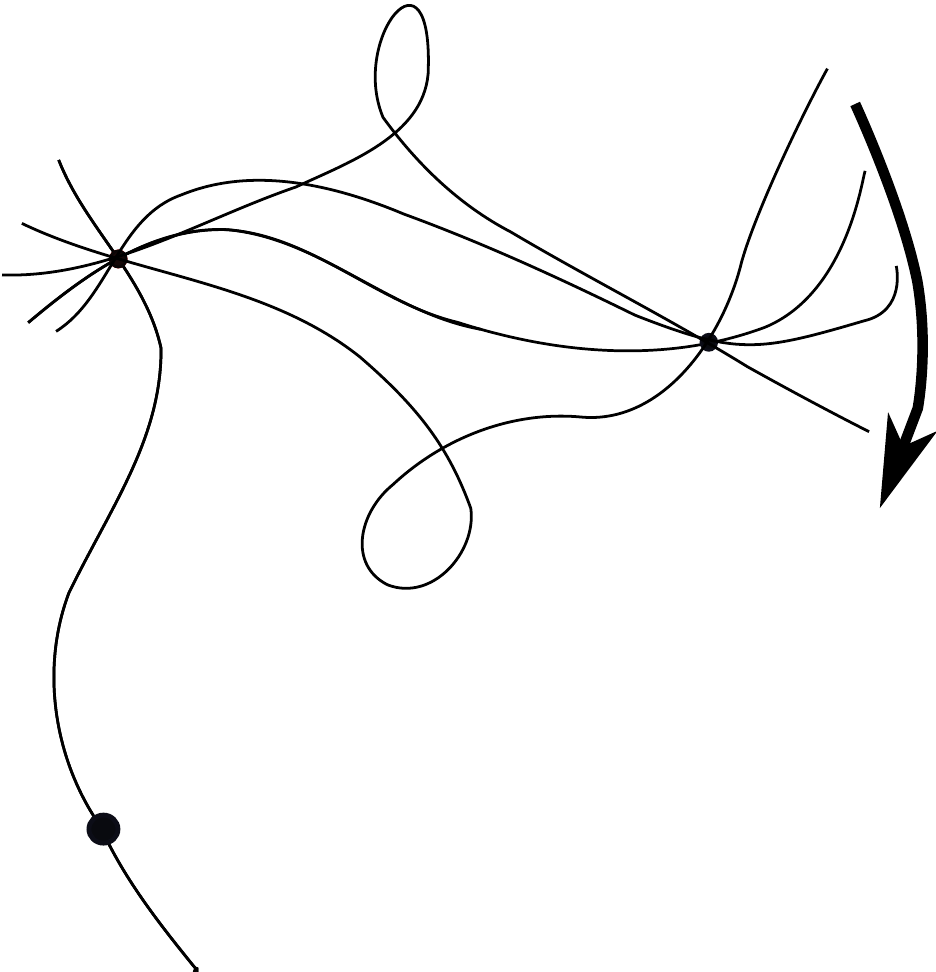}

\put (9,65){$y$}
\put (-35,18){$g(C)=g-1$}
\put (-5,35){$C$}
\put (60,35){A pencil of plane cubics}
\put (15,14){$q$}
\put (20,-5){\footnotesize{Figure 7: Test curve $E$ }}
\end{overpic}
\end{center}
\end{figure}

This is a standard test curve and it is well-known \cite{HarrisMorrison} that $E\cdot \lambda=1,E\cdot\delta_0=12,E\cdot\delta_{g-1}=-1$, giving
\begin{equation*}
E\cdot R=c_\lambda+12c_0-c_{g-1}
\end{equation*}
To find the intersection directly we observe that for any such solution either the limit $p$ is in the $C$-aspect or it is not. If the limit point $p$ lies on $C$ we have
\begin{equation*}
gp+q+(g-k-5)y+\sum_{j=1}^k q_j\sim K_C
\end{equation*} 
for some $k=1,...,g-3$. But any such solution would contradict our assumption that $q$ and $y$ are general. If $p$ does not lie on $C$ then we have
\begin{equation*}
q+(2g-k-5)y+\sum_{j=1}^k q_j\sim K_C
\end{equation*} 
for some $k=1,...,g-3$. Again we have a contradiction for any $k$ with the assumption that $q$ and $y$ are general points. Hence $E\cdot R=0$ and by our previous test curve results we see
\begin{equation*}
c_\lambda=\frac{g(3g^3-3g+2)}{2}
\end{equation*}
which agrees with our Porteous' formula result.

\bibliographystyle{plain}
\bibliography{base}
\end{document}